\documentclass{article}
\usepackage [latin1]{inputenc}
\usepackage{longtable,booktabs}
\usepackage{graphicx}
\usepackage{amsmath}
\usepackage{mathrsfs}
 \usepackage{cite}
\usepackage{color}
 \usepackage{indentfirst}
\begin{document}
\newcommand{\wzr}[1]{\textcolor{blue}{#1}}
\newcommand*{\doublerightarrow}[2]{\mathrel{
  \settowidth{\@tempdima}{$\scriptstyle#1$}
  \settowidth{\@tempdimb}{$\scriptstyle#2$}
  \ifdim\@tempdimb>\@tempdima \@tempdima=\@tempdimb\fi
  \mathop{\vcenter{
    \offinterlineskip\ialign{\hbox to\dimexpr\@tempdima+1em{##}\cr
    \rightarrowfill\cr\noalign{\kern.5ex}
    \rightarrowfill\cr}}}\limits^{\!#1}_{\!#2}}}
\baselineskip 0.25in
\title{\Large\bf MP and MT properties of fuzzy inference with aggregation function}
\author{{Dechao Li}\thanks{Email:\ dch1831@163.com},\qquad {Mengying He}\\
  {\small School of Information and Engineering, Zhejiang Ocean University,
  Zhoushan,
  316000, China}}
\date{}
\maketitle
\begin{center}
\begin{minipage}{140mm}
\begin{picture}(1,1)
 \line(1,0){400}
\end{picture}

\centerline{\bf Abstract} \vskip 3mm {\qquad
 As the two basic fuzzy inference models, fuzzy modus ponens (FMP) and fuzzy modus tollens (FMT) have the important application in artificial intelligence. In order to solve FMP and FMT problems, Zadeh proposed a compositional rule of inference (CRI) method. This paper aims mainly to investigate the validity of $A$-compositional rule of inference (ACRI) method, as a generalized CRI method based on aggregation functions, from a logical view and an interpolative view, respectively. Specifically, the modus ponens (MP) and modus tollens (MT) properties of ACRI method  are discussed in detail. It is shown that the aggregation functions to implement FMP and FMT problems provide more generality than the t-norms, uninorms and overlap functions as well-known the laws of $T$-conditionality, $U$-conditionality and $O$-conditionality, respectively. Moreover, two examples are also given to illustrate our theoretical results. Especially, Example 6.2 shows that the output $B'$ in FMP (FMT) problem is close to $B$ ($D^C$) with our proposed inference method when the fuzzy input and the antecedent of fuzzy rule are near (the fuzzy input near with the negation of the seccedent in fuzzy rule).}
 \vskip 2mm\noindent{\bf Key words}: Fuzzy implication; Aggregation function; ACRI method; MP and MT properties

\begin{picture}(1,1)
 \line(1,0){400}
\end{picture}
\end{minipage}
\end{center}
\vskip 4mm
\section{Introduction}
\subsection{Motivation of this paper}
Fuzzy reasoning has been successfully utilized in fuzzy control, artificial intelligence, data mining, image processing, decision-making, classification, prediction and so on \cite{Hullermeier, Kerre,Ma,Paul,Wang}. Since the exploration of fuzzy reasoning mechanism can help to improve the ability of human-like reasoning and thinking, it becomes an important topic to model the fuzzy reasoning in the applied science fields such as artificial intelligence.  As two basic inference models of fuzzy reasoning, FMP and FMT can draw some meaningful conclusions from some imprecise or vague premises using fuzzy inference method. In order to obtain the conclusions of FMP and FMT problems, the CRI method was proposed by Zadeh in 1973\cite{Zadeh1}. After, the inference
mechanism of CRI method was studied intensively by many researchers. And various t-norms and fuzzy implications in the CRI method
 were compared and investigated\cite{Baets,Bouchon,Driankov,Jenei,Mizumoto1,Morsi,Ruan}. Recently, Li and Zeng also extended the CRI method to the $A$-compositional rule of inference (ACRI) based on aggregation function\cite{Liz}. Unlike the CRI method, Pedrycz considered another inference method with the Bandler-Kohout subproduct (BKS) composition\cite{Pedrycz}. Moreover, other fuzzy reasoning inference methods such as the similarity-based approximate reasoning (SBR), triple implication principle (TIP) and quintuple implication principle (QIP) were also provided to solve the FMP and FMT problems\cite{Luo,Mizumoto,Pedrycz,Pei,Raha,Turksen,Wang1,Zhou}. It is worthy to mention that the CRI method is simpler than the other methods in computation. Therefore, it not only becomes one of the most
fundamental inference methods but also has been widely utilized in practice.

It is well known that the CRI method bases on the cylindrical extension and projection of fuzzy relations\cite{Zadeh1}. And
the t-norms and fuzzy implications have been employed to interpret the words ``and" and ``if...then..." in the CRI method, respectively. In order to meet the actual needs, some well-known fuzzy implications, which mainly include R-, $(S,N)$-, QL-, $f$-, $g$-, probabilistic-, probabilistic S-implications and $T$-power implications, were used to interpret the word ``if...then..." by researchers. Generally speaking, fuzzy implications can be generated as follows: i. by the binary functions on [0, 1], such as R-, $(S, N)$-, QL- and probabilistic implications \cite{Baczynski,Dimuro,Dimuro1,Grzegorzewski}; ii. by the unary functions on [0, 1], for instance, $f$- and $g$-implications \cite{Yager}. Considering the aggregation functions play a vital role in fuzzy logic, there also exist some fuzzy implications generated by aggregation functions such as the residual implications derived from aggregation functions and $(A,N)$-implications\cite{Ou,Pradera}.

However, the word ``and" cannot be precisely modeled by a single t-norm under all circumstances owing to the vagueness of natural language and the imprecision of people's thinking. Indeed, it does not accord with any t-norm now and then\cite{Zimmer1}. Moreover, the solution of CRI method $B'(y)(\equiv1)$ is unless or misleading in a singleton rule MISO (multi-input-single-output) fuzzy system when the t-norm is chosen to interpret the word ``and" while the fuzzy implication $I$ is employed to translate the fuzzy rule\cite{Li1,Li2}. The solution is rooted in which $T(a,0)\equiv0$ and $I(0,b)\equiv1$ always hold for a t-norm $T$ and a fuzzy implication $I$, respectively. This triggers us to seek some suitable operators correspond to the word ``and". Indeed, the operator extended from t-norm can flexibly effectively handle the uncertainties of words which cannot be handled by t-norms\cite{Hudec, Wu}. Moreover, considering the associativity or commutativity of operators to model the word ``and" is not necessarily required in decision making and classification problems, aggregation functions become increasingly concerned to substitute for the t-norms in the actual decision making and classification\cite{Bustince,Fodor}. Indeed, aggregation functions, as the generalization of t-norms and t-conorms, have been applied extensively in fuzzy logic, decision making and classification problems\cite{Bustince,Deveci, Dimuro,Fodor,Helbin,Liz,Mas,Mas1,Pradera,Pradera1,Pradera2}. Therefore, we need to model the word ``and" in the CRI method by an aggregation function. And then investigate the ACRI method for FMP and FMT problems to meet the actual demand.

We know that various tautologies can be regarded as inference rules to make deductive reasoning in classical logic. Modus ponens (MP) and Modus tollens (MT) are two most commonly inference rules. As the generalization of MP and MT in classical logic, FMP and FMT are used for making uncertain reasoning. With the intuition of human being, it is well-known that the conclusion $B'$ is close to the consequent of promise 1 when promise 2 is close to the antecedent of promise 1 in FMP problem. Especially, the conclusion is equal to the consequent of promise 1 if promise 2 is the antecedent of promise 1. Similarly, the farther are between promise 2 and the consequent of promise 1, the farther are between conclusion and the antecedent of promise 1. Especially, the conclusion is equal to the negation of consequent of promise 1 if promise 2 is the negation of antecedent of promise 1. Therefore, we need some basic requirements
for fuzzy reasoning to assess the effectiveness of inference methods for FMP and FMT problems\cite{Baldwin,Fukami,Magrez}.
\subsection{Literatures review}
Thereinto, Baldwin and Pilsworth represented four axioms from the logical view of fuzzy reasoning, including the classical MP property (see Section 2). On the other hand, Fukami et al. represented other four axioms from the interpolative view of fuzzy reasoning, including the classical MT property (see Section 2). Moreover, notice that the classical MP and MT properties have been studied for some fuzzy implications in the CRI method\cite{Alsina,Dimuro,Gera,Liw,Mas,Mas1,Mas2,Trillas}
For example, Alsina and Trillas investigated the $(T,T_1)$-condition of $(S,N)$-implications\cite{Alsina}. Dimuro et al. studied the $O$-conditionality of fuzzy implications based on the overlap and grouping functions\cite{Dimuro}.  Li et al. considered the Modus ponens property of $T$-power implications\cite{Liw}. Mas et al. discussed the $(U, N)$-
implications satisfying Modus Ponens\cite{Mas1}. Mas et al. also investigated the Modus Ponens and Modus Tollens
 of discrete implications\cite{Mas2}. Trillas provided the MPT-implication functions\cite{Trillas}. They are summarized in Table 1.\vspace{-2mm}
 $$\mbox{\bf{\small Table\ 1 \ Summary of refernene for MP and MT properties}}$$\vspace{-10mm}
\begin{center}
\begin{tabular}{ll}
 \toprule[1pt]
References &Contributions\vspace{1mm}\\
\midrule[0.75pt]
Baldwin and Pilsworth\cite{Baldwin}, 1980& Measure the validity of CRI method with MT\\
Fukami et al.\cite{Fukami}, 1980& Measure the validity of CRI method with MP and MT\\
Alsina and Trillas\cite{Alsina}, 2003 & Generalized MP for $(S,N)$-implications\\
Trillas\cite{Trillas}, 2004&MP and MT for the same implicaiton\\
Gera\cite{Gera}, 2007& MP and MT for special fuzzy sets\\
Mas et al. \cite{Mas2}, 2008&MP and MT for discrete fuzzy implication\\
Mas et al.\cite{Mas1}, 2018 &MP for $(U, N)$-implications\\
Mas et al. \cite{Mas}, 2019&MP for uninorm\\
Dimuro et al.\cite{Dimuro}, 2019&MP for overlap and grouping functions\\
Li et al. \cite{Liw}, 2022&MP for $T$-power implication\\
   \bottomrule[1pt]
\end{tabular}
\end{center}\vspace{1mm}
 Different from all methods above, we will employ the aggregation function to interpret the word ``and" in the FMP and FMT problems. And then investigate the validity of CRI method for well-known fuzzy implications from not only a logical view but also an interpolative view. To the best of our knowledge, there exists no works to research the MP and MT properties of ACRI method for more well-known fuzzy implications. This inspires us to investigate the validity of ACRI method with more well-known fuzzy implications from not only a logical view but also an interpolative view. Considering the ACRI method does not satisfy the axiom (A8) (see Theorem 4.2 in \cite{Liz}), our motivation of this paper is mainly to study the MP and MT properties of ACRI method with well-known fuzzy implications using the axioms (A4) and (A5). The specific discussions will be shown in Sections 4 and 5. Moreover, for the convenience of reading, some explanations of acronyms and symbols are provided in Tables 2 and 3.
 $$\mbox{\bf{\small Table\ 2 \ List of acronyms used in this paper }}\vspace{-3mm}$$
\begin{center}
\begin{tabular}{ll}
 \toprule[1pt]
Acronyms &Explanations\vspace{1mm}\\
\midrule[0.75pt]
FMP & Fuzzy modus ponens\\
 FMT &Fuzzy modus tollens \\
 CRI &Compositional rule of inference\\
  ACRI &$A$-compositional rule of inference\\
   MISO &Multi-input-single-output\\
   MP &Modus ponens\\
    MT &Modus tollens\\
LNC & Law of noncontradiction\\
LEM &Law of excluded middle\\
LB& Left boundary condition\\
RB& Right boundary condition\\
NP & Left neutrality property\\
IP &Identity principle\\
EP &Exchange principle\\
CP(N) &Law of contraposition with a fuzzy negation $N$\\
OP &Ordering property\\
LIA &Law of importation with $A$\\
DAC &Dual of $A$-conditionality with respect to $N$\\
AC &$A$-conditionality\\
RP&Residual propoerty\\
   \bottomrule[1pt]
\end{tabular}
\end{center}\vspace{1mm}
$$\mbox{\bf{\small Table\ 3 \ List of symbols used in this paper }}\vspace{-3mm}$$
\begin{center}
\begin{tabular}{ll}
 \toprule[1pt]
Symbols &Explanations\vspace{1mm}\\
\midrule[0.75pt]
$N_\bot$ & The smallest fuzzy negation\\
 $N_\top$ &The greatest fuzzy negation\\
 $N_c$ &The standard fuzzy negation\\
  $N_I$ &The natural negation of a fuzzy implication $I$\\
  $f_\varphi$ &The $\varphi$-conjugate of $f$\\
    $A$ &Aggregation function\\
$T$ &t-norm\\
$S$&t-conorm\\
$D_\top$ &The greatest disjunctor\\
$D_\bot$ &The smallest disjunctor\\
$S_{LK}$ & {\L}ukasiewicz t-conorm\\
$M_{\lambda,f}$ &Weighted quasi-arithmetic mean\\
$I$ &Fuzzy implication\\
$I_A$ &R-implication generated by $A$\\
$I_{A,N}$ &$(A,N)$-implication generated by $A$ and $N$\\
$I_{A_1,A_2}$&QL-operation generated by $A_1$ and $A_2$\\
$I_f$ & $f$-implication with an $f$-generator\\
$I_g$ & $g$-implication with a $g$-generator\\
$I^T$ &$T$-power implication\\
$I_C$ &Probabilistic implication\\
$\tilde{I}_C$ &Probabilistic S-implication\\
$D^C$ &The complement of $D$\\
$N_A$ & The natural negation of $A$\\
$C$&Copula\\
  $I^{lc}_{I,N}$ &$N$-lower-contrapositivisation of $I$\vspace{1mm}\\
   $I^{uc}_{I,N}$&$N$-upper-contrapositivisation of $I$\vspace{1mm}\\
   $d(D,D')$&Distance between $D$ and $D'$\\
   \bottomrule[1pt]
\end{tabular}
\end{center}\vspace{1mm}
 \subsection{Contributions of this paper}
Based on the argument above, we mainly examine the MP and MT properties of the ACRI method in this paper. Therefore, this paper firstly studies the properties of fuzzy implication and aggregation function in the ACRI method satisfying (A4) or (A5). And then the aggregation function is sought for well-known fuzzy implications in the ACRI method such that they satisfy (A5). Finally, we show the conditions of well-known fuzzy implications satisfying (A4) in the ACRI method and give an ACRI method involved in fuzzy implications and aggregation functions which satisfy (A4) or (A5). In short, this paper aims to:

(1) Investigate the properties of aggregation functions and fuzzy implications which satisfy (A4) and (A5) in the ACRI method, respectively.

(2) Construct the aggregation function for well-known fuzzy implications such that the ACRI method satisfies (A5).

(3) Show the conditions for the well-known fuzzy implications satisfying (A4) with a strong negation in the ACRI method.

(4) Propose an ACRI method with fuzzy implications and aggregation functions satisfying (A4) or (A5) to make inference better in the actual application.

The composition of this paper is as follows. Some basic concepts and definitions utilized in this paper is recalled in Section 2. Section 3 investigates the properties of fuzzy implications and aggregation functions satisfying (A4) or (A5) in the ACRI method. Section 4 constructs an aggregation function for well-known fuzzy implications mainly including R-, $(A,N)$-, QL-, $f$-, $g$-, probabilistic-, probabilistic S-implications and $T$-power implications  such that the ACRI method satisfies (A5), respectively.
In Section 5, the conditions for the well-known fuzzy implications satisfying (A4) with a strong negation in the ACRI method are investigated. Section 6 proposes a fuzzy inference method employed the fuzzy implication satisfying (A4) or (A5) and presents two examples to illustrate our theoretical results.
\section{Preliminary}
\subsection{Fuzzy negation, aggregation function and fuzzy implication}
In this section we will recall some main concepts used in this paper.
\\{\bf Definition 2.1}\cite{Lowen} A fuzzy negation function $N$ is a mapping $N: [0,1]\rightarrow [0,1]$ fulfilling

(N1) $N(0)=1, \ N(1)=0$;

(N2) $N(x)\geq N(y)\ \textmd{if}\ x\leq y,\ \forall\ x, y\in [0,1]$.

 Further, a fuzzy negation $N$ is strict if it satisfies the following properties:

(N3) $N$ is continuous;

(N4) $N(x)> N(y)\ \textmd{if}\ x<y$.

A fuzzy negation $N$ is strong if it is involutive, i.e.,

(N5)  $N(N(x))=x,
\forall\ x\in [0,1]$.
 \\{\bf Examples 2.2}\cite{Lowen}
 \begin{itemize}
   \item The following are the smallest and greatest fuzzy negations, respectively:

\qquad $N_\bot(x) = \left\{\begin{array}{ll}
               1 & x=0\\
               0& \textmd{otherwise}
             \end{array}\right.$ and\quad
$N_\top(x) = \left\{\begin{array}{ll}
               0 & x=1\\
               1& \textmd{otherwise}
             \end{array}\right.$.
\item The standard fuzzy negation $N_c(x)=1-x$.

   \item The natural negation of a fuzzy implication $I$ denoted by $N_I(x)=I(x,0)$ (See Definition 2.11).
 \end{itemize}

Let $\varphi$ be an automorphism on [0,1] i.e. an increasing bijection on [0,1] and $f$ be a binary function on [0,1]. The function $f_\varphi(x,y)=\varphi^{-1}(f(\varphi(x),\varphi(y))$ is called as the $\varphi$-conjugate of $f$.
\\{\bf Theorem 2.3}\cite{Klir} The fuzzy negation $N$ is strong if and only if there exists an automorphism $\varphi$ on [0,1] such that $N=(N_c)_\varphi$.
\\{\bf Definition 2.4} \cite{Grabisch} A binary aggregation function $A$ is a mapping $A:[0,1]^2\rightarrow [0,1]$  which meets the following conditions:

(A1) Boundary conditions: $A(0, 0)=0$ and $A(1, 1)= 1$;

(A2) Non-decreasing in each variable.

 Obviously, the $\varphi$-conjugate of aggregation function $A$, denoted by $A_\varphi$, is again an aggregation function.
\\{\bf Definition 2.5}\cite{Grabisch} For a binary aggregation function $A$, $e$ is referred as a left (right) neutral element if  $A(e, x)=x\ (A(x,e)=x)$ for all $x\in [0,1]$. Further, $e$ is called a neutral element if $A(e, x)=A(x,e)=x$.
\\{\bf Definition 2.6}\cite{Grabisch} A binary aggregation function $A$ is said to be

i. a conjunctor if $A(1, 0) =A(0, 1) =0$,

ii. a disjunctor if $A(1, 0) =A(0, 1) =1$,

iii. a semi-copula if 1 is a neutral element,

iv. associative if $A(x,A(y,z))=A(A(x,y),z)$,

v. commutative if $A(x,y)=A(y,x)$,

vi. a uninorm if it is associative, commutative and $e\in (0,1)$ is a neutral element,

vii.  a t-norm if it is an associative and commutative semi-copula,

ix. a t-conorm if it is the $N$-dual of a t-norm, that is, $S(x,y)=N^{-1}(T(N(x),N(y))$, where $T$ is a t-norm and $S$ a t-conorm.
\\{\bf Definition 2.7}\cite{Klement} The aggregation function $A_1$ is not greater than $A_2$ denoted  by $A_1\leq A_2$ if $A_1(x,y)\leq A_2(x,y)$ for any $x,y\in [0,1]$.
\\{\bf Definition 2.8} \cite{Pradera} Let $A$ be an aggregation function.

i. $A$ has zero divisors if there exist $x, y\in (0, 1]$ such that  $A(x,y) = 0$,

ii. $A$ has one divisors if there exist $x, y \in [0, 1)$ such that  $A(x,y) =1$.
\\{\bf Example 2.9}\cite{Klement,Grabisch,Pradera2} Some distinguished disjunctors are shown as follows:

\begin{itemize}
  \item The greatest disjunctor, $D_\top(x, y)=\left\{\begin{array}{ll}
                                               0 & x=y=0\\
                                               1&\textmd{otherwise}
                                             \end{array}\right.$;
  \item The smallest disjunctor, $D_\bot(x, y)=\left\{\begin{array}{ll}
                                               1 & x=1\ \textmd{or}\ y=1\\
                                               0&\textmd{otherwise}
                                             \end{array}\right.$;
  \item t-conorms, such as the {\L}ukasiewicz t-conorm $S_{LK}(x,y)=(x+y)\wedge 1$;
  \item $A(x,y) = f^{-1}((f(x\vee y)-f(N(x\wedge y)))\vee 0)$, where $N$ is a strong negation and $f:[0, 1]\rightarrow [0, +\infty]$ is a continuous and strictly decreasing with $f(1)=0$;

 \item Weighted quasi-arithmetic mean (WQAM), $M_{\lambda,f}(x, y)= f^{-1}((1-\lambda)f(x)+\lambda f(y))$, where $f:[0, 1]\rightarrow [-\infty, +\infty]$ is continuous strictly monotone with $f(1)=\pm\infty$ and $\lambda\in(0,1)$;
\item Group functions;
 \item $A_{e,g}(x,y) = g^{(-1)} (g(x)+ g(y))$,
where $g:[0, 1]\rightarrow [-\infty, \infty]$ is a continuous strictly monotone function with $g(e)=0$, $g(1)=\pm\infty$ and $g^{(-1)}$ its pseudo inverse of $g$ defined as\vspace{2mm} $g^{(-1)}(x)=\left\{\begin{array}{ll}
                                                                                    g^{-1}(x) & x\leq g(1) \\
                                                                                    1 & \textmd{otherwise}
                                                                                  \end{array}\right.$.
\end{itemize}

Considering the non-commutativity of aggregation function, we extend the non-contradiction principle in Ref.\cite{ Pradera2} as follows.
\\{\bf Definition 2.10}\cite{Pradera2} A binary aggregation function $A$ fulfills the law of non-contradiction (LNC) with respect to a fuzzy negation $N$ if
$$A(N(x),x)=A(x,N(x))=0,\ \forall x\in [0,1].\eqno(\textmd{LNC})$$

 By duality, $A$ fulfills the law
of excluded middle (LEM) with $N$ if
$$A(N(x),x)=A(x,N(x))=1, \ \forall x\in [0,1].\eqno(\textmd{LEM})$$
Evidently, (LNC) and (LEM) become the well-known law of non-contradiction and law
of excluded middle if $A$ is a t-norm and t-conorm, respectively. Moreover, $A$ is  a conjunctor (disjunctor) if it fulfills (LNC) ((LEM)).
\\{\bf Definition 2.11}\cite{Baczynski1}  A function  $I: [0, 1]^2\rightarrow[0, 1]$ is called a fuzzy implication if

(I1) Non-increasing in the first variable, i.e. $I(x,z)\geq I(y, z)$ if $x\leq y$;

(I2) Non-decreasing in the second variable, i.e. $I(x,y)\leq I(x, z)$ if $y\leq z$;

(I3) $I(0,0)=1$;

(I4) $I(1,1)=1$;

(I5) $I(1,0)=0$.

According to Definition 2.11, a fuzzy implication fulfills the following properties:

(LB) Left boundary condition, $I(0, y)= 1, \forall\ y\in[0, 1]$;

(RB) Right boundary condition, $I(x, 1)=1, \forall\ x\in[0, 1]$.
\\{\bf Definition 2.12} \cite{Baczynski1} A fuzzy implication  $I: [0, 1]^2\rightarrow[0, 1]$ satisfies for any $x,y,z\in [0,1]$:

(NP) Left neutrality property if $I(1, y)=y$;

(IP) Identity principle if $I(x, x) = 1$;

(EP) Exchange principle if $I(x, I(y, z)) = I(y, I(x, z))$;

(CP(N)) Law of contraposition with a fuzzy negation $N$ if $I(x, y) = I(N(y),N(x))$;

(OP) Ordering property if $I(x, y)=1\Longleftrightarrow x\leq y$.
\\{\bf Definition 2.13}\cite{Pradera} Let $A$ be an aggregation function and $I$ a fuzzy implication. $I$ satisfies the law of importation with $A$ (LIA) if
$$I(A(x,y),z)=I(x,I(y,z)), \ \ \forall x, y,z \in [0, 1].\eqno(\textmd{LIA})$$
Obviously, (LIA) become the well-known law of importation if $A$ is a t-norm. It is also proved that (LIA) can be used to construct a hierarchical fuzzy systems (see Ref.\cite{Lig}).
\\{\bf Definition 2.14}\cite{Baczynski1} Let $I$ and $J$ be two fuzzy implications. Defining $(I\vee J)(x,y)=\max(I(x,y),$ $J(x,y))$, $(I\wedge J)(x,y)=\min(I(x,y),J(x,y))$ and $I_\varphi(x,y)=\varphi^{-1}(I(\varphi(x),\varphi(y)))$,  $I\vee J$, $I\wedge J$ and $I_\varphi$ are fuzzy implications, too.
\\{\bf Definition 2.15} \cite{Ou}  An R-implication generated by the aggregation function $A$ is a function $I_A: [0, 1]^2 \rightarrow [0, 1]$ defined as $$I_A(x,y) = \sup\{ t \in [0, 1]\mid A(x, t) \leq y\}.$$
{\bf Definition 2.16}\cite{Pradera}  A function $I_{A,N}:[0,1]^2\rightarrow[0,1]$ is called an $(A,N)$-implication generated by a disjunctor $A$ and a fuzzy negation $N$ defined as $$I_{A,N}(x,y)=A(N(x), y).$$
 $I_{A,N}$ is an $A$-implication if $N$ is the standard negation.  Moreover, an $(A,N)$-implication generated by a t-conorm and a strong negation is called a strong implication (S-implication).
\\{\bf Theorem 2.17}\cite{Pradera} $I$ is a fuzzy implication if and only if there exists a disjunctor $A$ such that $I(x,y)=A(1-x,y)$.
\\{\bf Definition 2.18}\cite{Pradera} A QL-operation $I_{A_1,A_2}$ generated by two aggregation functions $A_1$, $A_2$ and a fuzzy negation $N$ is defined as
$$I_{A_1,A_2}(x, y) = A_1(N(x), A_2(x, y)) ,\quad\quad x, y \in [0, 1].$$
 Especially, it is called a QL-implication if it fulfills (I1) and (I3)-(I5).
 \\{\bf Definition 2.19}\cite{Yager} Let a mapping $f:[0,1]\rightarrow [0,+\infty]$ be strict decreasing and continuous with
$f(1)=0$. An $f$-implication with an $f$-generator is defined by
$I_{f}(x,y)=f^{-1}(xf(y))$ with the understanding $0\times \infty=0$.
\\{\bf Definition 2.20}\cite{Yager} Suppose that $g:[0,1]\rightarrow [0,+\infty]$ is a strict increasing and continuous mapping with
$g(0)=0$. A $g$-implication with a $g$-generator is defined by
$I_{g}(x,y)=g^{(-1)}\left(\frac{g(y)}{x}\right)$ with the understanding $0\times \infty=\infty$.
\\{\bf Definition 2.21}\cite{Massanet} A $T$-power implication is defined by  $I^T(x, y) = \vee\{r\in [0, 1] | y^{(r)}
_T \geq x\}$ for all $x, y\in[0, 1]$, where $T$ is a continuous t-norm $T$.
\\{\bf Lemma 2.22}\cite{Massanet} Let $I^T$ be a $T$-power implication. Then\vspace{1mm}

i. $I^{T_\textmd{M}}(x,y)=\left\{\begin{array}{cc}
                                                                                                 1 & x\leq y\\
                                                                                                 0& x>y
                                                                                               \end{array}\right.$;

ii. $I^T(x, y) =\left\{\begin{array}{ll}
                                                                                                 1 & x\leq y\\
                                                                                                 \frac{t(x)}{t(y)}& x>y
                                                                                               \end{array}\right.$, where $T$ is Archimedean with additive generator $t$.\vspace{1mm}
\\{\bf Definition 2.23}\cite{Grzegorzewski}  A probabilistic implication $I_C$ generated by a copula $C$ is defined by\vspace{2mm} $I_C(x,y)=\left\{\begin{array}{ll}
                                                                                   \frac{C(x,y)}{x} & x>0 \vspace{1mm}\\
                                                                                    1 & \textmd{otherwise}
                                                                                  \end{array}\right.$
if it satisfies (I1).\vspace{1mm}
\\{\bf Definition 2.24}\cite{Grzegorzewski} A probabilistic S-implication $\widetilde{I}_C$ generated by a copula $C$ is defined as
$\tilde{I}_C(x,y)=C(x,y)-x+1$.
\\{\bf Definition 2.25}\cite{Zadeh} Let $D$ be a fuzzy set on the universe $U$. We say that $D$ is a normal fuzzy set if there exist some $x_0\in U$ such that $D(x_0)=1$.
\subsection{Fuzzy inference and axioms systems}
We know that FMP and FMT problems can be expressed as follows:
$$\textmd{Premise\ 1:}\ \textmd{IF}\ x\ \textmd{is }D\ \textmd{THEN}\ y\ \textmd{is}\ B
\qquad\qquad \textmd{Premise\ 1:}\ \textmd{IF}\ x\ \textmd{is }D\ \textmd{THEN}\ y\ \textmd{is}\ B\vspace{-2mm}$$
$$\textmd{Premise\ 2:}\ \ x\ \textmd{is }D'\qquad\quad\qquad\qquad\qquad\qquad \textmd{Premise\ 2:}\qquad\quad\qquad\qquad\ \ y\ \textmd{is }B'\vspace{-3mm}$$
\begin{picture}(1,2)
 \line(1,0){170}\qquad
 \line(1,0){170}\qquad
\end{picture}\qquad\vspace{-3mm}
$$\textmd{Conclusion:}\qquad\quad\qquad\qquad y\ \textmd{is}\ B',\qquad\qquad \textmd{Conclusion:}\ x\ \textmd{is}\ D',\qquad\qquad\qquad\qquad\qquad$$
where $D$ and $D'$ are fuzzy sets on the universe $U$ while $B$ and $B'$ are fuzzy sets on the universe $V$.

In order to calculate the conclusion $B'$ in FMP problem, we will substitute an aggregation function $A$
for $T_M$, and constrain the fuzzy relation to the fuzzy implication in the CRI method in this paper. Then, the following ACRI method is considered:
$$B'_{\textmd{ACRI}}(y)=\bigvee_{x\in U}A(D'(x),I(D(x), B(y))),$$
 where  $I$ is a fuzzy implication  and $A$ is an aggregation function.

There are many different axiom systems as the basic requirements for fuzzy inference methods. The following are the most well-known axiom systems:
\begin{itemize}
  \item The axioms presented by Baldwin and Pilsworth\cite{Baldwin}.

(A1)\quad $B\subseteq B'$;

(A2)\quad $B'\subseteq B''$ if $D'\subseteq D''$;

(A3)\quad $B'=V$ if $D'=D^C$ with $D^C$ being the complement of $D$;

(A4)\quad $B'=D^C$ if $D'=B^C$(classical modus tollens property);

\item The axioms presented by Fukami et al.\cite{Fukami}.

(A5)\quad $B'=B$ if $D'=D$ (classical modus ponens property);

(A6)\quad $B'=very\ B$ if $D'= very\ D$;

(A7)\quad $B'= more\ or\ less\ B$ if $D'= more\ or\ less\ D$;

(A8)\quad $B'=B^C$ if $D'= D^C$.
\end{itemize}

In this paper we follow not only the logic view but also interpolative view of the ACRI method. Therefore, our first aim is to investigate the properties of aggregation functions and fuzzy implications satisfying (A4) and (A5) in the ACRI method, respectively.
\section{Properties of aggregation functions and fuzzy implications satisfying (A4) or (A5)}
  Considering the normal fuzzy sets play an important role in the ACRI method, we always assume the fuzzy sets involved in the ACRI method are normal in the remainder of this paper. Let us firstly investigate the properties of aggregation function and fuzzy negation in the ACRI method fulfilling (A4) for a  fuzzy implication. For convenience, we say that the fuzzy implication $I$ satisfies (A4) with an aggregation function $A$ and a fuzzy negation $N$ if the ACRI method fulfills (A4) from now on.
\\{\bf Proposition 3.1} Let the fuzzy implication $I$ satisfy (A4) with an aggregation function $A$ and a fuzzy negation $N$. Then, $A$ is a conjunctor.
\\{\bf Proof.} Since $I$ satisfies (A4) with $A$ and $N$, we have $N(D(x))=\mathop{\bigvee}\limits_{y\in V}A(N(B(y)),$ $I(D(x),B(y)))$ for all $x\in U$ and $y\in V$. Especially, $A(0,I(1,1))=A(0,1)=0$ holds for $x_0\in U$ and $y_0\in V$ such that $D(x_0)=1$ and $B(y_0)=1$. Similarly, we can obtain $A(1,0)=0$. Thus, $A$ is a conjunctor.
\\{\bf Proposition 3.2} Let the fuzzy implication $I$ fulfill (NP) and $A$ be a conjunctor having a right neutral element 1. If $I$ satisfies (A4) with $A$ and a fuzzy negation $N$, then $N_I\leq N\leq N_A$, where $N_A(x)=\sup\{t\in[0,1]|A(t,x)=0\}$ is the natural negation of $A$.
 \\{\bf Proof.} Since $I$ satisfies (A4) with $A$ and $N$, $N(D(x))\geq A(N(0),I(D(x),0))=N_I(D(x))$ holds for any $x\in U$. We can similarly obtain $0=N(1)=A(N(B(y)),$ $I(1,B(y)))=A(N(B(y)),B(y))$. This implies that $N\leq N_A$.
 \\{\bf Proposition 3.3} Let the fuzzy implication $I$ fulfill (IP) and the aggregation function $A$ have a right neutral element 1. Then, $I$
satisfies (A4) with $A$ and $N$ if and only if $A(N(B(y)),I(D(x),$ $B(y)))\leq N(D(x))$ holds for all $x\in U$ and $y\in V$.
\\{\bf Proof.} $(\Longrightarrow)$ Obviously.

$(\Longleftarrow)$ For any $x\in U$ and $y\in V$, we assume that $A(N(B(y)),I(D(x),B(y)))\leq N(D(x))$ holds. This implies that $\mathop{\bigvee}\limits_{y\in V}A(N(B(y)),I(D(x),B(y)))\leq N(D(x))$. On the other hand,\vspace{1mm} $\mathop{\bigvee}\limits_{y\in V}A(N(B(y)),$ $I(D(x),B(y)))\geq A(N(D(x)),I(D(x),D(x)))$. Since $I$ fulfills (IP) and 1 is a\vspace{1mm} right neutral element of $A$,  we have $A(N(D(x)),I(D(x),$ $D(x)))=A(N(D(x)),1)=N(D(x))$. Therefore, $I$
satisfies (A4) with $A$ and $N$.
\\{\bf Proposition 3.4} Let the aggregation function $A$ have a right neutral element 1. Then, the fuzzy implication $I$
satisfies (A4) with $A$ and $N_I$ if and only if $A(N_I(B(y)),I(D(x),B(y)))\leq N_I(D(x))$ holds for all $x\in U$ and $y\in V$.
\\{\bf Proof.} $(\Longrightarrow)$ Obviously.

$(\Longleftarrow)$ For any $x\in U$, we have $\mathop{\bigvee}\limits_{y\in V}A(N_I(B(y)),I(D(x),B(y)))\geq A(N_I(0),$ $I(D(x),0))=A(1,N_I(D(x)))=N_I(D(x))$.

Proposition 3.4 shows that the inequality  $A(N_I(B(y)),I(D(x),B(y)))\leq N_I(D(x))$ can replace the axiom (A4) in this case. Therefore, we say that the fuzzy implication $I$ satisfies the dual of $A$-conditionality with respect to $N$ if the above inequality holds.  It is not difficult to see that $B(y)$ and $D(x)$ variate on [0,1]. And $N_I$ is a fuzzy negation (see Example 2.2). This means that they can be denoted by some variables (such as $a$ and $b$) in [0,1].  And then this inequality is shorten as
$$A(N(b),I(a,b))\leq N(a), \ \forall \ a,b\in[0,1].\eqno(\textmd{DAC})$$
{\bf Proposition 3.5} Let the fuzzy implication $I$ meet the condition $I(a,0)=0$ for all $a\in (0,1]$ and $N=N_\bot$. Then, $I$ satisfies (DAC) with $N_\bot$ and any conjunctor $A$.
\\{\bf Proof.} It is sufficient to verify that $A(N_\bot(b),I(a,b))=0$ holds for any $a\neq 0$. Let us consider the following two cases.

i. $b\neq 0$. This implies that $A(N_\bot(b),I(a,b))=A(0,I(a,b))=0$ holds for any $a\in(0,1]$.

 ii.  $b=0$. In this case, we have $A(1,I(a,0))= A(1,0)=0$ for any $a\in(0,1]$.
\\{\bf Proposition 3.6} Let $N$ be an injective fuzzy negation and the aggregation function $A$ have a right neutral element 1. If there exist $1>a>b>0$ such that $I(a,b)=1$, then the fuzzy implication $I$ does not satisfy (A4) with $A$.
\\{\bf Proof.} Let $1>a_0>b_0>0$ and $I(a_0,b_0)=1$. This implies that $A(N(b_0),I(a_0,$ $b_0))=A(N(b_0),1)=N(b_0)>N(a_0)$ holds.
\\{\bf Proposition 3.7} Let the fuzzy implication $I$ satisfy (A4) with an aggregation function $A$  and a fuzzy negation $N$. If $N\neq N_{\perp}$ and $I(1,b)>0$ holds for any $b>0$, then $A$ has zero divisors.  Especially, $A$ satisfies (LNC) if $I$ fulfills (NP).
\\{\bf Proof.} By (A4), we have $A(N(B(y)),I(1, B(y)))=0$ for all $y\in V$. This implies that $A$ has zero divisors.

Next, we study the properties of aggregation function and fuzzy implication in the ACRI method satisfying (A5). Similarly, we say that the fuzzy implication $I$ satisfies (A5) with an aggregation function $A$ when the ACRI method  satisfying (A5) for convenience.
\\{\bf Proposition 3.8} Let the aggregation function $A$ have a right neutral element 1. If the fuzzy implication $I$ satisfies (A5) with $A$, then $I(1,b)\leq b$ holds for all $b\in[0,1]$.
\\{\bf Proof.} Since $I$ satisfies (A5), we have $B(y)=\mathop{\bigvee}\limits_{x\in U}A(D(x),I(D(x),B(y)))$ for any $y\in V$. This\vspace{1mm} implies that $B(y)\geq A(1,I(1,B(y))=I(1,B(y))$ for any $y\in V$. Since $B$ is normal,  $I(1,b)\leq b$ holds for all $b\in [0,1]$.
\\{\bf Proposition 3.9} Let the fuzzy implication $I$ fulfill (NP) and the aggregation function $A$ have a left neutral element 1. $I$ satisfies (A5) with $A$ if and only if $A(D(x),I(D(x),B(y)))\leq B(y)$ holds for all $x\in U$ and $y\in V$.
\\{\bf Proof.} $(\Longrightarrow)$ Obviously.

$(\Longleftarrow)$ Suppose that $A(D(x),I(D(x),B(y)))\leq B(y)$ holds for any $x\in U$ and $y\in V$. We immediately have $\mathop{\bigvee}\limits_{x\in U}A(D(x),I(D(x),B(y)))\leq B(y)$ for any $y\in V$. On the other\vspace{1mm} hand, there exists $x_0\in U$ such that $D(x_0)=1$. Therefore, $\mathop{\bigvee}\limits_{x\in U}A(D(x),I(D(x),B(y)))\geq$\vspace{1mm} $ A(D(x_0),I(D(x_0),B(y)))=A(1,I(1,B(y)))=I(1,B(y))=B(y)$ holds for any $y\in V$.
\\{\bf Remark 1.} i. Similar to Ref.\cite{Baczynski1,Dimuro,Liw,Mas1}, we say that the fuzzy implication $I$ satisfies $A$-conditionality if the above inequality holds. It is easy to find that $B(y)$ and $D(x)$ variate on [0,1]. Thus, they can be denoted by some variables (such as $a$ and $b$) in [0,1]. And then, this inequality is shorten as
$$A(a,I(a,b))\leq b, \ \forall \ a,b\in[0,1].\eqno(\textmd{AC})$$

ii. Let the fuzzy implication $I$ fulfill (IP) and the aggregation function $A$ have a right neutral element 1. Then, $I$ satisfies (A5) with $A$ if they fulfill (AC).

iii. Obviously, (AC) implies (DAC) if the fuzzy implication $I$ satisfies (CP(N)). Further, (DAC) implies (AC) if the fuzzy implication $I$ fulfills (CP(N)) and the fuzzy negation $N$ is continuous.
Owing to the relationship between (AC) and (DAC), the following results can be parallelly obtained and similarly proved.
\\{\bf Proposition 3.10} If the fuzzy implication $I$ satisfies (A5) with an aggregation function $A$, then $A$ is a conjunctor.
\\{\bf Proposition 3.11} Let the fuzzy implication $I$ satisfy (A5) with an aggregation function $A$. Then, $A$ has zero divisors if $N_I>N_{\perp}$.
\\{\bf Proposition 3.12} Let the aggregation function $A$ have a right neutral element 1. If there exist $1>a>b>0$ such that $I(a,b)=1$, then the fuzzy implication $I$ does not satisfy (A5) with $A$.
\\{\bf Remark 2.} i. Similar results can be found in Ref.\cite{Baczynski1} if the aggregation function $A$ becomes a t-norm (see Propositions 7.4.2 and 7.4.3).

 ii. According to Propositions 3.6 and 3.12, $I(a,b)=1$ implies $a\leq b$ if the fuzzy implication $I$ satisfies (A4) or (A5) with $A$ having a right neutral element 1. However, we cannot ensure that the fuzzy implication $I$ satisfies (OP) in this case as shown in the following examples.
\\{\bf Example 3.13} Let $I_f$ be an $f$-implication. Obviously, $I_f$ does not fulfill (OP). Define an aggregation function $A_f$ as $A_f(x,y)=\left\{\begin{array}{ll}
                                                                                                       f^{-1}(\frac{f(y)}{x})& x\neq 0\\
                                                                                                       0& x=0
                                                                                                    \end{array}\right.$. Then, $(I_f,A_f)$ satisfies (A5).\vspace{1mm}
\\{\bf Example 3.14} Let $I$ be an R-implication generated by a nilpotent t-norm. Obviously, $I$ fulfills (OP). However, $I$ does not satisfies
(A5) with $A$ and $N_\top$ unless $A$ is the smallest conjunctor.
\\{\bf Proposition 3.15} Let the fuzzy implication $I$ fulfill (OP) and $A$ be a commutative aggregation function. If $I$ satisfies (LIA) with $A$, then $I$ satisfies (AC) with $A$.
\\{\bf Proof.} Since $I$ satisfies (LIA) with $A$ and $I$ fulfills (OP), we obtain the fact that $I$ satisfies (EP) by Lemma 3.1 in Ref.\cite{Lig}. Then, we have  $1=I(I(a,b),I(a,b))=I(a,I(I(a,b),b))=I(A(a,I(a,b)),b)$ for all $a,b\in[0,1]$. Therefore, $A(a,I(a,b))\leq b$ holds for all $a,b\in[0,1]$.
\\{\bf Proposition 3.16} Let $I$ and $J$ be two fuzzy implications fulfilling (NP) and $J\leq I$. If $I$ satisfies (A5) with $A$ having a left neutral element 1, then $J$ satisfies (A5) with $A$, too.
\\{\bf Proof.} By Proposition 3.8, it is sufficient to verify that $J$ satisfies (AC) with $A$. Indeed, the monotonicity of $A$ implies that
$b\geq A(a,I(a,b))\geq A(a,J(a,b))$.
\\{\bf Remark 3.} We can similarly obtain the fact that $J$ satisfies (DAC) with $A$ and $N$ if $I$ satisfies (DAC) with $A$ and $N$. However, we cannot ensure whether $J$ satisfies (A4) with $A$ and $N$ when $J$ satisfies (A4) with $A$ and $N$.
\\{\bf Proposition 3.17} Let the fuzzy implication $I$ fulfill (NP) and the aggregation function $A$ have  a left neutral element 1. We have

i. if $A'$ is a conjunctor having a left neutral element 1 such that $A'\leq A$ and $I$ satisfies (A4) with $A$ and $N_I$, then $I$ satisfies (A4) with $A'$ and $N_I$, too;

ii. if $A'$ is a conjunctor having a left neutral element 1 such that $A'\leq A$ and $I$ satisfies (A5) with $A$, then $I$ fulfills (A5) with $A'$, too.
\\{\bf Proof.} We only verify $I$ satisfies (A4) with $A'$ and $N$. The other case can be verified similarly. By Proposition 3.3, it is sufficient to verify that $I$ fulfills (DAC) with $A'$. Indeed,
$b\geq A(a,I(a,b))\geq A'(a,I(a,b))$ holds.
\\{\bf Remark 4.} Baczynski presented the similar results if the aggregation function $A$ becomes a t-norm in Ref.\cite{Baczynski1} (see Proposition 7.4.3).

We know that the $\varphi$-conjugate of the fuzzy implication $I$ still is a fuzzy implication. Thus, it is interesting to investigate the  satisfaction of $I_\varphi$ for (A4) and (A5) when $I$ fulfills respectively (A4) and (A5).
\\{\bf Proposition 3.18} i. If the fuzzy implication $I$ satisfies (A4) with an aggregation function $A$ and a fuzzy negation $N$, then $I_\varphi$ fulfills (A4) with $A_\varphi$ and $N_\varphi$, where $A_\varphi$, $I_\varphi$ and $N_\varphi$ are the $\varphi$-conjugate of $A$, $I$ and $N$, respectively;

ii. If the fuzzy implication $I$ satisfies (A5) with an aggregation function $A$, then $I_\varphi$ fulfills (A5) with $A_\varphi$.
\\{\bf Proof.} We only verify $I_\varphi$ satisfies (A5) with $A_\varphi$. The other can be verified similarly. For any $y\in V$, we have $\mathop{\bigvee}\limits_{x\in U}A_\varphi(D(x), I_\varphi(D(x),B(y)))=\mathop{\bigvee}\limits_{x\in U}\varphi^{-1}(A(\varphi(D(x)), $\vspace{1mm} $\varphi(\varphi^{-1}(I(\varphi(D(x)),\varphi(B(y)))))))=\mathop{\bigvee}\limits_{x\in U}\varphi^{-1}(A(\varphi(D(x)), I(\varphi(D(x)),\varphi(B(y)))))=\varphi^{-1}(\mathop{\bigvee}\limits_{x\in U}A(\varphi(D(x)), I(\varphi(D(x)),$ $\varphi(B(y)))))=\varphi^{-1}(\varphi(B(y)))=B(y)$.

Similarly, $I\vee J$ and $I\wedge J$ still are two  fuzzy implications if $I$ and $J$ are fuzzy implications. Thus, it is necessary to study the  satisfaction of (A4) and (A5) for $I\vee J$ and $I\wedge J$ when $I$ and $J$ fulfill (A4) and (A5) respectively.
\\{\bf Proposition 3.19} Let $I$ and $J$ be two fuzzy implications. We have

i. If both $I$ and $J$ satisfy (A4) with an aggregation function $A$ and a fuzzy negation $N$, then $I\vee J$ and $I\wedge J$ fulfill respectively (A4) with $A$ and $N$, too.

ii. If both $I$ and $J$ satisfy (A5) with an aggregation function $A$, then $I\vee J$ and $I\wedge J$ satisfy respectively (A5) with $A$, too.
\\{\bf Proof.} Obviously.
\section{MP property of the ACRI method with well-known fuzzy implications}
As the argument above, the aggregation functions and fuzzy implications are respectively employed to interpret the word ``and" and ``if...then..." in the ACRI method. Thus, this section will discuss the MP property of the ACRI method for some well-known fuzzy implications including R-, $f$-, $g$-, $T$-power, $(A,N)$-, QL-, probabilistic- and probabilistic S-implications. That is, we will construct aggregation functions such that these well-known fuzzy implications satisfy (A5) with them, respectively. Then, we firstly have the following result.
\\{\bf Lemma 4.1} Let the fuzzy implication $I$ fulfill the condition $I(1,b)<1$ for any $b\in[0,1)$. Then, there exists an aggregation function $A_I$ defined as $ A_I(a,b)=\inf\{c\in[0,1]|I(a,c)\geq b\}$ such that $(I,A_I)$ satisfies (AC).
\\{\bf Proof.} According to Lemma 3.3 in \cite{Liz}, $A_I$ is an aggregation function. Further, we have $b\geq\inf\{c\in[0,1]|I(a,c)\geq I(a,b)\}= A_I(a,I(a,b))$.
\\{\bf Corollary 4.2} Let the fuzzy implication $I$ fulfill (NP). Then $(I,A_I)$ satisfies (A5).
\\{\bf Proof.} By Proposition 3.9 and Lemma 4.1, it is sufficient to verify that $A_I$ has a left neutral element 1. Indeed, $A_I(1,b)=\inf\{c|I(1,c)\geq b\}=\inf\{c|c\geq b\}=b$ holds for any $b\in [0,1]$.

Next, we construct an aggregation function for the R-implication involved in the ACRI method such that they satisfy (A5).  As the R-implication $I_A$ generated by a left-continuous aggregation function with respect to second variable fulfills the residuation property (RP), we then have the following statement.
\\{\bf Theorem 4.3} Let $A$ and $A'$ be two conjunctors having a left neutral element 1. If $A$ is left-continuous with respect to second variable, then the R-implication $I_A$ satisfies (A5) with $A'$ if and only if $A'\leq A$.
 \\{\bf Proof.} $(\Longleftarrow)$ Since 1 is a left neutral element of $A$, the R-implication $I_A$ fulfills (NP). According to Proposition 3.9, it is sufficient to verify that $(I_A,A')$ satisfies (AC).  The left-continuity of  $A$ with respect to second variable implies that $A$ and $I_A$ satisfy the residuation property (RP) by Lemma 3.1 in \cite{Liz}, i.e. $A(a,c)\leq b\Longleftrightarrow c\leq I_A(a,b)$ holds for any $a,b,c\in[0,1]$. Especially, we have $I_A(a,b)\leq I_A(a,b)\Longleftrightarrow A(a,I_A(a,b))\leq b$ by (RP). This means that $b\geq A(a,I_A(a,b))\geq A'(a,I_A(a,b))$ holds for all $a,b\in [0,1]$.

$(\Longrightarrow)$ Suppose that $I_A$ satisfies (A5) with $A'$. We have $A'(a,I_A(a,b))\leq b$ for any $a,b\in[0,1]$. The residuation property of $I_A$ with $A$ implies that $b\leq I_A(a,A(a,b))$ for any $a,b\in [0,1]$. Therefore, $A(a,b)\geq A'(a,I_A(a,A(a,b)))\geq A'(a,b)$ holds for any $a,b\in[0,1]$.

As it can be easily seen, it is a strict condition to require the aggregation function $A$ being left-continuous with respect to second variable in  R-implication $I_A$. Therefore, we next focus on the the R-implication generated by a non-left-continuous aggregation function with respect to second variable. It is not difficult to see that the R-implication $I_A$ generated by an aggregation function $A$ having a right neutral element 1 satisfies $I_A(a,b)=1$ if $a\leq b$. This means that the R-implication $I_A$ fulfills (OP) when $I_A$ satisfies (A5) with $A'$ according to Proposition 3.12. In order to investigate the satisfaction of (A5) for the R-implication $I_A$, it is sufficient to study this case where the R-implication $I_A$ fulfills (OP) and $A$ has a right neutral element 1. We then have the following statement.
\\{\bf Lemma 4.4} Let $I_A$ be an R-implication generated by an aggregation function $A$ having a right neutral element 1. Then, $I_A$ satisfies (OP) if and only if $A$ is border continuous (that is, $A$ is continuous at the border of $[0,1]^2$).
\\{\bf Proof}. This can be proved similarly to Proposition 2.5.9 in \cite{Baczynski1}.

 Inspired by the idea of \cite{Jayaram}, we will extend a border continuous aggregation function to a left-continuous one.
 \\{\bf Definition 4.5} Let $A$ be a border continuous aggregation function. A function $A^*$ on $[0,1]^2$ is defined as follows:
 \begin{equation}A^*(x,y)=\left\{\begin{array}{ll}
  \sup\{A(u,v)|u<x,v<y\}& x,y\in(0,1) \vspace{1mm}\\
  A(x,y) & \textmd{otherwise}
  \end{array}\right..\end{equation}
Obviously, $A^*$ is a left-continuous aggregation function with respect
to second variable and $A^*\leq A$.
 \\{\bf Theorem 4.6} Let $I_A$ be an R-implication generated by an aggregation function $A$ having a neutral element 1 and the conjunctor $A'$ have a left neutral element 1. If $I_A$ fulfills (OP) and (EP), then $I_A$ satisfies (A5) with $A'$ if and only if $A'\leq A^*$ holds, where $A^*$ is defined as Eq.(1).
 \\{\bf Proof.}  We firstly assert that $I_A=I_{A^*}$. Obviously, $I_A\leq I_{A^*}$ holds. We further have $I_A(a,b)=I_{A^*}(a,b)=1$ if $a\leq b$. Thus, it needs to consider the case $1>a>b>0$. On the contrary, suppose that there exist $1>a_0>b_0>0$ such that $c=I_A(a_0,b_0)<I_{A^*}(a_0,b_0)=c'$. By intermediate value theorem, there exists $c_0$ such that $c<c_0<c'$. Since $A^*$ is a left-continuous aggregation function with respect to second variable, $c_0\leq I_{A^*}(a_0,b_0)=c'$ implies that $A^*(a_0,c')\leq b_0$ holds. This means $b_0\geq A^*(a_0,c')=\sup\{A(u,v)|u<a_0,v<c'\}\geq \sup\{A(u,c_0)|u<a_0\}
 =A(a_0^-,c_0)$. Therefore, we have $c_0\leq I_A(a_0^-,b_0)$. Then, $1=I_A(c_0, I_A(a_0^-,b_0))=I_A(a_0^-, I_A(c_0,b_0))$ holds by (OP) and (EP). Again, we obtain $a_0^-\leq I_A(c_0,b_0)$ according to (OP). On the other hand, we have $1>I_A(c_0, I_A(a_0,b_0))=I_A(a_0, I_A(c_0,b_0))$ by (OP) and (EP). Thus, $a_0> I_A(c_0,b_0)$ holds. This implies that $I_A(c_0,b_0)=a_0$. However, $I_A(c_0,c)=I_A(c_0, I_A(a_0,b_0))=I_A(a_0, I_A(c_0,b_0))=I_A(a_0,a_0)=1$ by (OP) and (EP). Again, we can obtain $c_0\leq c$. This is a contradiction.

 For the $f$-, $g$- and  $T$-power implications, it is easy to construct the aggregation functions  according to Lemma 4.1. Therefore, we immediately have the following results by Corollary 4.2.
\\{\bf Theorem 4.7} Let $I_f$ be an $f$-implication with an $f$-generator. Then, $I_f$ satisfies (A5) with an\vspace{1mm} aggregation function if and only if $A\leq A_{I_f}$, where $A_{I_f}(x,y)=\left\{\begin{array}{ll}
f^{-1}(\frac{f(y)}{x})& x\neq0\\
0& x=0
\end{array}\right.$ for any $x,y\in[0,1]$.\vspace{1mm}
\\{\bf Proof.} $(\Longleftarrow)$ Since $I_f$ fulfills (NP), $I_f$ satisfies (A5) with $A$ if $A\leq A_{I_f}$ by Corollary 4.2.

$(\Longrightarrow)$  $A_{I_f}(a,I_f(a,b))=\inf\{c\in [0,1]|I_f(a,c)\geq I_f(a,b)\}=b$ holds for any $a\in(0,1]$ and $b\in[0,1]$. Since $(I_f,A)$ satisfies (A5), we have $A((a,I_f(a,b))\leq b=A_{I_f}(a,I_f(a,b))$ for any $a\in(0,1]$ and $b\in[0,1]$. The continuity of $I_f$ implies that $A(a,b)\leq A_{I_f}(a,b)$ holds for any $a\in(0,1]$ and $b\in[0,1]$. Setting $a=0$, we have $A(0,b)=0=A_{I_f}(0,b)$. Thus,  $A\leq A_{I_f}$ holds.
\\{\bf Theorem 4.8} Let $I_g$ be a $g$-implication with a $g$-generator. Then $I_g$ fulfills (A5) with an\vspace{1mm} aggregation function $A$ if and only if $A\leq A_{I_g}$,\vspace{1mm} where
$A_{I_g}(x,y)=\left\{\begin{array}{ll}
g^{-1}(xg(y))& x\neq0\\
0& x=0
\end{array}\right.$ for any $x,y\in[0,1]$.\vspace{1mm}
\\{\bf Proof.} This can be proved similarly to Theorem 4.7.
 \\{\bf Theorem 4.9} Let $I^T$ be a $T$-power implication.

i. $I^{T_M}$ satisfies (A5) with an aggregation function $A$ if and only if $A$ has a right neutral element 1;

ii. If $T$ is an Archimedean t-norm with additive generator $t$, then $I^T$ satisfies (A5) with an\vspace{1mm} aggregation function $A$ if and only if $A\leq A_{I^T}$, where $A_{I^T}(x,y)=
\left\{\begin{array}{ll}\vspace{1mm}
t^{-1}(\frac{t(x)}{y})&y\neq 0\\
0&\textmd{otherwise}
\end{array}\right.$.\vspace{1mm}
\\{\bf Proof.} This can be proved similarly to Theorem 4.7.

However, it is difficult to obtain the aggregation function $A_I$ for $(A,N)$-, QL-, probabilistic and probabilistic S-implications by Lemma 4.1. So, we only consider some special $(A,N)$-, QL-, probabilistic and probabilistic S-implications in the last of this section. Let us firstly consider the cases when $(A, N)$-implication generated by the smallest disjunctor $D_{\bot}$.
 \\{\bf Proposition 4.10} Let $I_{D_{\bot},N}$ be an $(A, N)$-implication generated by the smallest disjunctor $D_{\bot}$. Then, $I_{D_{\bot},N}$ does not satisfy (A5) with any conjunctor.
 \\{\bf Proof.} The $(A, N)$-implication generated by the smallest disjunctor $D_{\bot}$  is the smallest fuzzy\vspace{1mm} implication  $I_{D_{\bot},N}(x,y)=\left\{\begin{array}{ll}\vspace{1mm}
1& x=0\ \textmd{or}\ y=1\\
0& \textmd{otherwise}
\end{array}\right.$. Let $A$ be a conjunctor. For any fuzzy set $B$\vspace{1mm} on $V$ such that $B(y)\neq 0$, we have $B(y)>\mathop{\bigvee}\limits_{x\in U}A'(D(x),I_{D_{\bot},N}(D(x),B(y)))=\mathop{\bigvee}\limits_{x\in U}A(D(x),0)=0$.\vspace{1mm}
\\{\bf Proposition 4.11}  Let $I_{D_{\top},N}$ be an $(A, N)$-implication generated by the greatest disjunctor $D_{\top}$.
 Then $I_{D_{\top},N}$ does not satisfy (A5) with any aggregation function.\vspace{1mm}
 \\{\bf Proof.} In this case, $I_{D_{\top},N}(x,y)=\left\{\begin{array}{ll}\vspace{1mm}
0& x=1,y=0\\
1& \textmd{otherwise}
\end{array}\right.$. Let $A$ be an aggregation function. For\vspace{1mm} any $y\in V$ such that $B(y)\in (0,1)$, we have $A(1,$\vspace{1mm} $I_{D_{\top},N}(1,B(y)))=1>B(y)$.

Next, we study the case where $(A, N)$-implication is generated by a continuous t-conorm with the ordinal sum structure. According to Corollary 5.12 in \cite{Klement}, there exist a uniquely determined index set $\Gamma$, a set of uniquely determined open pairwise
disjoint intervals $\{(a_\alpha, e_\alpha)\}_{\alpha\in \Gamma}$ of [0, 1] and a set of uniquely
determined Archimedean continuous t-conorms $(S_\alpha)_{\alpha\in \Gamma}$ such that $S$ can be rewritten as
$$S(x,y)=\left\{\begin{array}{ll}
a_\alpha+(e_\alpha-a_\alpha)S_\alpha\left(\frac{x-a_\alpha}{e_\alpha-a_\alpha}, \frac{y-a_\alpha}{e_\alpha-a_\alpha}\right)& x,y\in[a_\alpha, e_\alpha]\vspace{1mm}\\
x\vee y& \textmd{otherwise}\\
\end{array}\right..$$
Then, we have the following statement.
\\{\bf Theorem 4.12}  Let $I_{S,N}$ be an $(A, N)$-implication generated by a continuous t-conorm $S$ with the ordinal sum structure and a continuous fuzzy negation $N$. Then,
$I_{S,N}$ fulfills (A5) with an aggregation function $A$ if and only if $A\leq A_{I_{S,N}}$, where $ A_{I_{S,N}}$ is defined as
$$A_{I_{S,N}}(x,y)=\left\{\begin{array}{ll}
0& N(x)\geq y\\
a_\alpha+(e_\alpha-a_\alpha)f^{-1}_\alpha\left(f_\alpha\left(\frac{N(x)-a_\alpha}{e_\alpha-a_\alpha}\right)
-f_\alpha\left(\frac{y-a_\alpha}{e_\alpha-a_\alpha}\right)\right)& N(x)\in [a_\alpha, e_\alpha]\ \textmd{and}\ N(x)<y\\
y& \textmd{otherwise}\\
\end{array}\right.$$ and $f_\alpha$ is the continuous additive generator of Archimedean t-conorm $S_\alpha$.\vspace{1mm}
\\{\bf Proof.}
With some tedious calculations according to Lemma 4.1, we can obtain $A_{I_{S,N}}$  as follows

$A_{I_{S,N}}(x,y)=\left\{\begin{array}{ll}
0& N(x)\geq y\\
a_\alpha+(e_\alpha-a_\alpha)f^{-1}_\alpha\left(f_\alpha\left(\frac{N(x)-a_\alpha}{e_\alpha-a_\alpha}\right)
-f_\alpha\left(\frac{y-a_\alpha}{e_\alpha-a_\alpha}\right)\right)& N(x)\in [a_\alpha, e_\alpha]\ \textmd{and}\ N(x)<y\\
y& \textmd{otherwise}\\
\end{array}\right.$,\vspace{1mm} where $f_\alpha$ is the additive generator of continuous Archimedean t-conorm $S_\alpha$.

$(\Longleftarrow)$ Since $I_{S,N}$ fulfills (NP), $I_{S,N}$ satisfies (A5) with $A$ if $A\leq A_{I_{S,N}}$ by Corollary 4.2.

$(\Longrightarrow)$ Suppose that $I_{S,N}$ satisfies (A5) with an aggregation function $A$.  We consider the following three cases:

i. If $N(a)\geq b$. In this case, we have $A(a,b)\leq A(a,N(a))=0=A_{I_{S,N}}(a,b)$.

ii. If $N(a)\in [a_\alpha, e_\alpha]$ and $N(a)<b$. This case implies that $A_{I_{S,N}}(a,I_{S,N}(a,b))=b\geq A(a,I_{S,N}(a,b))$ holds. By the continuity of $I_{S,N}$, we have $A(a,b)\leq A_{I_{S,N}}(a,b)$.

iii. If $N(a)\notin [a_\alpha, e_\alpha]$ and $N(a)<b$. In this case, we have $I_{S,N}(a,b)=b$. Therefore,  $A(a,b)=A(a, I_{S,N}(a,b)) \leq b=A_{I_{S,N}}(a,b)$.
\\{\bf Remark 5.} We can similarly obtain the results about the $(A, N)$-implications generated by the dual of representable aggregation function, weighted quasi-arithmetic mean, group functions generated by continuous functions with a neutral element $e$, respectively. Here, the repetitious details are shown no longer.

For a QL-operation $I_{A_1,A_2}$, it is easy to see that $I_{A_1,A_2}$ satisfies (I3) and (I5) when $A_1$ is a disjunctor and $A_2$ is a conjunctor. Further, let the conjunctor $A_2$ have a left neutral element 1. Then, we have the fact that $A_1$ satisfies (LEM) if $I_{A_1,A_2}$ is a QL-implication. Therefore, we only consider the case where $I_{A_1,A_2}$ is obtained from a disjunctor $A_1$ having a left neutral element 0, a conjunctor $A_2$ having a neutral element 1 and a fuzzy negation $N$ in the remainder of this section.
\\{\bf Lemma 4.13} Let $I_{A_1,A_2}$ be a QL-implication generated by a disjunctor $A_1$ without one divisor, a conjunctor $A_2$ and a fuzzy negation $N$ mentioned above. Then, $I_{A_1,A_2}$ fulfills (A5) with any conjunctor having a left neutral element 1.
\\{\bf Proof.}  Since $A_1$ has not one divisor, $A_1$ fulfills (LEM) with $N$ if and only if $N=N_\top$.  Therefore, $I_{A_1,A_2}$ is an $(A,N)$-implication generated by $A_1$ and $N_\top$. We can easily verify that $I_{A_1,A_2}$ satisfies (A5) with conjunctor $A$ having a left neutral element 1.

Further, we have $I_{A_1,A_2}(x,y)=A_1(N(x),A_2(x,y))\leq A_1(N(x),A_2(1,y))=A_1(N(x),y)=I_{A_1,N}(x,y)$. This means that $I_{A_1,A_2}$ satisfies (A5) with the same conjunctor $A$ if $I_{A_1,N}$ fulfills (A5) with $A$. In a general way,  we can find
a disjunctor $A$ such that $I_{A_1,A_2}(x,y)=I_{A,N}(x,y)=A(1-x,y)$ according to Theorem 2.17. Thus, we have the following statement.
\\{\bf Theorem 4.14} Let $I_{A_1,A_2}$ be a QL-implication. Then, there exists always an $(A,N)$-implication $I_{A,N}$ such that $I_{A_1,A_2}$ satisfies (A5) with the same conjunctor $A'$ if and only if $I_{A,N}$ fulfills (A5) with $A'$.
\\{\bf Proof.} Obviously.
\\{\bf Theorem 4.15} Let $C$ be an Archimedean copula with additive generator $c$. Then,

 i. The probabilistic implication $I_C$ generated by $C$ satisfies (A5) with $A$ if and only if\vspace{1mm} $A\leq A_{I_C}$, where $A_{I_C}$ is\vspace{1mm} defined as $A_{I_C}(x,y)=$ $\left\{\begin{array}{ll}
0& x=0\ \textmd{or}\ y=0\\
c^{-1}(c(xy)-c(x))& \textmd{otherwise}
\end{array}\right.$ and $c$ is the additive generator of $C$;\vspace{1mm}

 ii. The probabilistic S-implication $\widetilde{I}_C$ generated by $C$ satisfies (A5) with $A$ if and only if\vspace{1mm} $A\leq A_{\widetilde{I}_C}$, where $A_{\widetilde{I}_C}(x,y)=$$\left\{\begin{array}{ll}
0& x+y\leq1\\
c^{-1}(c(x+y-1)-c(x))& \textmd{otherwise}
\end{array}\right.$.\vspace{1mm}
\\{\bf Proof.} The proof is similar to that of Theorem 4.12.
\section{MT property of the ACRI method with well-known fuzzy implications}
  In this section we will investigate the MT property of the ACRI method for well-known fuzzy implications. Considering the fuzzy negation plays a vital role when the ACRI method satisfies (A5), we firstly study the case where the fuzzy negations are $N_\bot$ and $N_\top$, respectively. Obviously, we have firstly the following statements.
 \\{\bf Lemma 5.1} Let $A$ be a conjunctor and the  fuzzy implication $I$ fulfill the condition $I(a,0)=0$. If $N_I=N_\bot$, then $I$ satisfies (A4) with $A$ and $N_\bot$.
  \\{\bf Proof.} By Lemma 3.5, $I$ satisfies (DAC) with $N_I$ and any conjunctor $A$. Therefore, it is sufficient to verify that the equals sign holds in (DAC). Indeed, $\mathop{\bigvee}\limits_{y\in V}
A(N_\bot(B(y)),I(D(x),B(y)))=N_I(D(x))$.
 \\{\bf Lemma 5.2} Let the fuzzy implication $I$ fulfill the condition $I(1,b)>0$ for some $b\in (0,1)$. Then, there does not exist any conjunctor such that $I$ satisfies (A4) with the greatest fuzzy negation $N_{\top}$.
 \\{\bf Proof.} This proof is similar to that of Lemma 5.1.

However, it is difficult to investigate the satisfaction of (A4) in the ACRI method with ordinary non-continuous fuzzy negation. We therefore only consider the case where $N$ is a strong fuzzy negation in the rest of this section. As pointed out in Remark 1, the (CP(N)) acts as a bridge between the MP and MT properties in the ACRI method. Thus, we assume that the fuzzy implication $I$ fulfills (NP) and the aggregation function $A'$ has a left neutral element 1. And then we study the (CP(N)) with a strong fuzzy negation $N$ for $(A,N)$-, R-, QL-, $f$-, $g$-, probabilistic-, probabilistic S- and $T$-power implications, respectively.
\\{\bf Proposition 5.3} Let $N$ be a strong fuzzy negation. Then, $(A,N)$-implication $I_{A,N}$ satisfies (CP(N)) if and only if $A$ is commutative.
\\{\bf Proof.} Obviously.
\\{\bf Lemma 5.4} Let $A$ be a commutative and left-continuous aggregation function with respect to second variable. Then, the R-implication $I_{A}$ satisfies (CP(N)) if and only if $I_A(x,y)=N_I(A(x,N_I(y)))$.
\\{\bf Proof.} $(\Longrightarrow)$ Suppose that  $I_{A}$ fulfills (CP(N)) with a strong fuzzy negation $N$. We have $N=N_I$ according to Corollary 1.5.7 in\cite{Baczynski1}. The left-continuity of $A$ with respect to second variable implies that $A$ and $I_A$ satisfy (RP) by Lemma 3.1 in \cite{Liz}. Therefore, we have $A(x,y)\leq z\Longleftrightarrow y\leq I_A(x,z)\Longleftrightarrow y\leq I_A(N_I(z),N_I(x))\Longleftrightarrow A(N_I(x),y)\leq N_I(z)\Longleftrightarrow z\leq N_I(A(N_I(x),y))$. This implies that $I_A(x,y)=\max\{z\in [0,1]|A(x,z)\leq y\}=\max\{z\in [0,1]|z\leq N_I(A(N_I(x),y))\}=N_I(A(N_I(x),y))$.

$(\Longleftarrow)$ Since $A$ is commutative and $N_I$ is strong, we can easily verify that $I_A(x,y)=N_I(A(x,N_I(y))$ satisfies (CP(N)).
\\{\bf Remark 6.} i. In this case, $I_A(x,y)=A_{N_I}(N_I(x),y)$ is an $(A,N)$-implication, where $A_{N_I}(x,y)=N(A(N_I(x),N_I(y)))$ is the dual of $A$ with $N_I$.

ii. Obviously, the aggregation function $A$ has zero divisors. And then $A(x,y)=0$ holds if and only if $x\leq N_I(y)$.

For a QL-implication $I_{A_1,A_2}$, the following equation holds if it fulfills (CP(N)):
$$A_1(N(x), A_2(x, y))=A_1(y, A_2(N(x),N(y))).$$

Evidently, the disjunctor $A_1$ satisfies (LEM) if the conjunctor $A_2$ has a right neutral element 1. However, it is still not easy to solve this equation. So, we only consider the case when $A_1$ is a continuous t-conorm and $A_2$ is a t-norm. And the corresponding result can be found in \cite{Fodor1}.
\\{\bf Lemma 5.5}\cite{Baczynski1}  $I_f$ fulfills (CP(N)) with a strong fuzzy negation if and only if $f(0)<\infty$.
\\{\bf Lemma 5.6}\cite{Baczynski1} $I_g$ does not satisfy (CP(N)) with any fuzzy negation.
\\{\bf Lemma 5.7}\cite{Baczynski} $I_C$ does not satisfy (CP(N)) with any fuzzy negation.
\\{\bf Lemma 5.8}\cite{Baczynski} $\tilde{I}_C$ satisfies (CP(N)) if and only if the equation
$C(x, y)= x + y-1 + C(1-y, 1-x)$ holds.
\\{\bf Lemma 5.9}\cite{Massanet2} i. If $T$ is the Minimum t-norm, then $I^T$ satisfies (CP(N)) if and only if $N$ is strictly decreasing;

ii. If $T$ is a strict t-norm with additive generator $t$, then $I^T$ satisfies (CP(N)) if and only if  $N$ is  a strong fuzzy  negation given by $N(x)=t^{-1}(\frac{k}{t(x)}) $ for  some positive constant $k$;

iii. If $T$ is  a non-strict Archimedean t-norm, then $I^T$ does not satisfies (CP(N)) with any fuzzy negation.
\\{\bf Theorem 5.10} Let the fuzzy implication $I$ be R-implication, $(A,N)$-implication, QL-implication, $f$-implication, probabilistic S-implication and $T$-power implications fulfilling (CP(N)), respectively. If $I$ satisfies (A5) with an aggregation function $A$, then $I$ fulfills (A4) with the same $A$.
\\{\bf Proof.} We only consider the case where $I$ is an R-implication. Other cases can be similarly proved. Since $I$ fulfills (CP(N)), we have $\mathop{\bigvee}\limits_{y\in V}
A(N(B(y)),I(D(x),$ $B(y)))=\mathop{\bigvee}\limits_{y\in V}
A(N(B(y)),$ $I(N(B(y)),N(D(x))))=N(D(x))$. In the last equality, we use the fact that $I$ satisfies (A5) with $A$.

However, it is not easy that the fuzzy implications fulfill (CP(N)). Therefore, it is worthy to consider the case when the fuzzy implications do not satisfy (CP(N)). We firstly modify the fuzzy implication by the method in \cite{Aguilo} such that it satisfies (CP(N)). Let $I$ be a fuzzy implication and $N$  a strong fuzzy negation. The $N$-lower (upper)-contrapositivisation of $I$, denoted as $I^{lc}_{I,N}$ ($I^{uc}_{I,N}$)  is defined as
$$I^{lc}_{I,N}(x,y)=\left\{\begin{array}{ll}
I(x,y)& y\geq N(x)\vspace{1mm} \\
I(N(y),N(x))& \textmd{otherwise}
\end{array}\right.,\quad I^{uc}_{I,N}(x,y)=\left\{\begin{array}{ll}
I(x,y)& y\leq N(x)\vspace{1mm} \\
I(N(y),N(x))& \textmd{otherwise}
\end{array}\right..$$
We can obtain the fact that $I^{lc}_{I,N}$ and $I^{uc}_{I,N}$ satisfy (CP(N)) (see Theorems 1 and 4 in \cite{Aguilo}). Further, the aggregation functions $A_{I^{lc}_{I,N}}$ and $A_{I^{uc}_{I,N}}$ can be constructed respectively according to Lemma 4.1.  We have immediately the following result.
\\{\bf Theorem 5.11} Let $I$ be a fuzzy implication and $N$ a strong fuzzy negation. $I^{lc}_{I,N}$ and $I^{uc}_{I,N}$ satisfy (DAC) with  $A_{I^{lc}_{I,N}}$ and $A_{I^{uc}_{I,N}}$, respectively.
\\{\bf Proof.} This proof comes from Lemma 4.1.

Obviously, $I^{lc}_{I,N}$ and $I^{uc}_{I,N}$ are not equal to $I$ unless $I$ fulfills (CP(N)). For these well-known fuzzy implications, we begin to investigate the satisfaction of (A4) in case where the aggregation function is a special conjunctor.  As Proposition 3.7 implies that $I$ does not satisfy (A4) with any conjunctor $A$ without zero divisor, we thus move on the next case when $A$ has zero divisors. Inspired by Ref.\cite{Fodor}, we further demand that the aggregation function $A$ meets the following conditions:

(C1) $A$ has a left neutral element 1;

(C2) $A(x,1)<1$ and $\phi(x)=A(x,1)$ is strictly continuous increasing;

(C3) $A(x,A(y,z))=A(y,A(x,z))$.

According to Theorem 5.1 in \cite{Fodor}, there exists a t-norm $T$ such that $A(x,y)=T(\phi(x),y)$. In this case, we can easily see that $A$ has zero divisors if $T$ has zero divisors. Therefore, we have the following statements.
\\{\bf Theorem 5.12} Let $A$ be a continuous conjunctor fulfilling (C1)-(C3). Then, the fuzzy implication $I$ satisfies (DAC) with $A$ and $N$ if and only if there exists an automorphism $\varphi$ on [0,1] such that $I(a,b)\leq \varphi^{-1}(\varphi(N(a))+1-\varphi(\phi(N(b))))\wedge1)$. Especially, $I$ satisfies (A4) with $A$ and $N_I$ if and only if there exists an automorphism $\varphi$ on [0,1] such that $I(a,b)\leq \varphi^{-1}(2-\varphi(a)-\varphi(\phi(\varphi^{-1}(1-\varphi(b))))\wedge1)$ holds for any $a,b\in[0,1]$.
\\{\bf Proof.} $(\Longleftarrow)$ Obviously.

$(\Longrightarrow)$ Suppose that $I$ satisfies (DAC) with $A$ and $N$. The continuity of $A$ implies that $I(a,b)\leq I_A(N(b),N(a))$ holds for any $a,b\in[0,1]$, where $I_A$ is the R-implication generated by $A$.  Thus, $I_A(N(b),N(a))=\bigvee\{z\in[0,1]|A(N(b),z)\leq N(a)\}=
\bigvee\{z\in[0,1]|T(\phi(N(b)),z)\leq N(a)\}=\bigvee\{z\in[0,1]|\varphi^{-1}((\varphi(\phi(N(b)))+\varphi(z)-1)\vee 0)\leq N(a)\}=\varphi^{-1}(\varphi(N(a))+1-\varphi(\phi(N(b))))\wedge1)$. Therefore, $I(a,b)\leq \varphi^{-1}(\varphi(N(a))+1-\varphi(\phi(N(b))))\wedge1)$.

Especially, $I$ satisfies (A4) with $A$ and $N_I$ if 1 is  a left  neutral element of $A$ by Proposition 3.4. This means that  $I$ satisfies (A4) with $A$ and $N_I$ if and only if there exists an automorphism $\varphi$ on [0,1] such that $I(a,b)\leq \varphi^{-1}(2-\varphi(a)-\varphi(\phi(\varphi^{-1}(1-\varphi(b))))\wedge1)$ holds for any $a,b\in[0,1]$.
\\{\bf Theorem 5.13} Let $I_{A,N}$ is an $(A,N)$-implication generated by a disjunctor $A$ with a right neutral element 0 and a strong negation $N$. If the continuous conjunctor $A'$ fulfills (C1)-(C3),  then $I_{A,N}$ satisfies (A4) with $A'$ and $N$ if and only if there exists an automorphism $\varphi$ on [0,1] such that $A(a,b)\leq \varphi^{-1}(\varphi((N(a))+1-\varphi(\phi(N(b))))\wedge1)$ holds for any $a,b\in[0,1]$. Especially, $I_{A,N}$ satisfies (A4) with $A'$ and $N_I$ if and only if there exists an automorphism $\varphi$ on [0,1] such that $A(a,b)\leq (S_{LK})_\varphi(a,N_\varphi(\phi(N(b))))$ holds for any $a,b\in[0,1]$, where $(S_{LK})_\varphi$ is the $\varphi$-conjugate of {\L}ukasiewicz t-conorm $S_{LK}$.
\\{\bf Proof.} $(\Longleftarrow)$ This can be verified directly.

$(\Longrightarrow)$ Assume that $I_{A,N}$ satisfies (A4) with $A'$ and $N$.  Since $A$ has a right neutral element 0, we have $N_{I_{A,N}}=N\leq N_{A'}=N_\varphi$. Further, there exists an automorphism $\varphi$ on [0,1] such that $(\varphi(\phi(N(b)))+\varphi(A(a,b))-1)\vee0\leq\varphi(a)$.
This implies that $\varphi(A(a,b))\leq (1-\varphi(\phi(N(b)))+\varphi(a))\wedge 1$. Then, we have $A(a,b)\leq \varphi^{-1}(1-\varphi(\phi(N(b)))+\varphi(a))\wedge 1)$. Therefore, $A(a,b)\leq (S_{LK})_\varphi(a,N_\varphi(\phi(N(b))))$ holds for any $a,b\in[0,1]$, where $(S_{LK})_\varphi$ is the $\varphi$-conjugate of {\L}ukasiewicz t-conorm $S_{LK}$.
\\{\bf Theorem 5.14} Let $I_T$ be an R-implication generated by a continuous t-norm $T$. If the continuous conjunctor $A$ fulfills (C1)-(C3), then $I_T$ satisfies (DAC) with $A$ and $N$ if and only if there exist an automorphism $\varphi$ on [0,1] and some continuous additive generators $t_\alpha$ such that the following inequality holds for any $a,b\in[a_\alpha,e_\alpha]$:
$$t_\alpha\left(\frac{b-a_\alpha}{e_\alpha-a_\alpha}\right)
-t_\alpha\left(\frac{a-a_\alpha}{e_\alpha-a_\alpha}\right)\leq t_\alpha\left( \frac{\varphi^{-1}(\varphi(N(a))+1-\varphi(\phi(N(b))))\wedge1)-a_\alpha}{e_\alpha-a_\alpha}\right).$$
Especially, $I_T$ satisfies (A4) with $A$ and $N_I$ if and only if there exist an automorphism $\varphi$ on [0,1]  and some continuous additive generators $t_\alpha$ such that the following inequality holds for any $a,b\in[a_\alpha,e_\alpha]$:
$$t_\alpha\left(\frac{b-a_\alpha}{e_\alpha-a_\alpha}\right)
-t_\alpha\left(\frac{a-a_\alpha}{e_\alpha-a_\alpha}\right)\leq t_\alpha\left(  \frac{\varphi^{-1}(2-\varphi(a)-\varphi(\phi(\varphi^{-1}(1-\varphi(b))))\wedge1)-a_\alpha}{e_\alpha-a_\alpha}\right).$$
{\bf Proof.} $(\Longleftarrow)$ This can be verified directly.

$(\Longrightarrow)$ Let $I_T$ satisfy (A4) with $A$ and $N$. According to Theorem 5.12, we have $I_T(a,b)\leq \varphi^{-1}(\varphi(N(a))+1-\varphi(\phi(N(b))))\wedge1)$ for all $a,b\in[0,1]$. By Theorem 2.5.24 in \cite{Baczynski1}, $I_T$ can be rewritten as follows:
$$I_T(x,y)=\left\{\begin{array}{ll}
                    1 & x\leq y\\
                    a_\alpha+(e_\alpha-a_\alpha)I_{T_{\alpha}}\left(\frac{x-a_\alpha}{e_\alpha-a_\alpha},
                    \frac{y-a_\alpha}{e_\alpha-a_\alpha}\right)&x,y\in[a_\alpha,e_\alpha]\\
y&\textmd{otherwise}\\
                  \end{array}
\right.,$$
where $I_{T_\alpha}$ is an R-implication generated by the Archimedean t-norm $T_\alpha$. Obviously, it is sufficient to study the case when $a,b\in[a_\alpha,e_\alpha]$. The continuity of $T_\alpha$ implies that there exists a continuous additive generator $t_\alpha$ such that $T_\alpha(x,y)=t_\alpha^{-1}((t_\alpha(x)+t_\alpha(y))\wedge t(0))$. We therefore\vspace{1mm} have $t_\alpha\left(\frac{b-a_\alpha}{e_\alpha-a_\alpha}\right)
-t_\alpha\left(\frac{a-a_\alpha}{e_\alpha-a_\alpha}\right)\leq t_\alpha\left( \frac{\varphi^{-1}(\varphi(N(a))+1-\varphi(\phi(N(b))))\wedge1)-a_\alpha}{e_\alpha-a_\alpha}\right)$.\vspace{1mm}

Similarly, we have the following statements for QL-, $f$-, $g$-, probabilistic, probabilistic S- and $T$-power implications.
\\{\bf Theorem 5.15} Let $I_{A_1,A_2}$ be a QL-implication. If the continuous conjunctor $A$ fulfills (C1)-(C3), then $I_{A_1,A_2}$ satisfies (DAC) with $A$ and $N$ if and only if there exists an automorphism $\varphi$ on [0,1] such that $A_1(a, A_2(N(a), b))\leq \varphi^{-1}(\varphi(a)+1-\varphi(\phi(N(b))))\wedge1)$ holds for any $a,b\in[0,1]$.
Especially, $I_{A_1,A_2}$ satisfies (A4) with $A$ and $N_I$ if and only if there exists an automorphism $\varphi$ on [0,1] such that
$A_1(a, A_2(N(a), b))\leq \varphi^{-1}(1+\varphi(a)-\varphi(\phi(\varphi^{-1}(1-\varphi(b))))\wedge1)$  holds for any $a,b\in[0,1]$.
\\{\bf Proof.} $(\Longleftarrow)$ Obviously.

$(\Longrightarrow)$ We assume that $I_{A_1,A_2}$ satisfies (DAC) with $A$ and $N$. The continuity of $A$ implies that $I_{A_1,A_2}(a,b)\leq I_A(N(b),N(a))$ holds for any $a,b\in[0,1]$, where $I_A$ is the R-implication generated by $A$.  By Theorem 5.12, we have $A_1(a, A_2(N(a), b))\leq \varphi^{-1}(\varphi(a)+1-\varphi(\phi(N(b))))\wedge1)$.

Especially, $I_{A_1,A_2}$ satisfies (A4) with $A$ and $N_I$ if 1 is  a left  neutral element of $A$ by Proposition 3.4. This means that  $I_{A_1,A_2}$ satisfies (A4) with $A$ and $N_I$ if and only if there exists an automorphism $\varphi$ on [0,1] such that $A_1(a, A_2(N(a), b))\leq \varphi^{-1}(\varphi(a)+1-\varphi(\phi(N(b))))\wedge1)$ holds for any $a,b\in[0,1]$.
Therefore, $I_{A_1,A_2}$ satisfies (A4) with $A$ and $N_I$ if and only if there exists an automorphism $\varphi$ on [0,1] such that
$A_1(a, A_2(N(a), b))\leq \varphi^{-1}(1+\varphi(a)-\varphi(\phi(\varphi^{-1}(1-\varphi(b))))\wedge1)$  holds for any $a,b\in[0,1]$.
\\{\bf Theorem 5.16} Let $I_f$ be an $f$-implication. If the continuous conjunctor $A$ fulfills (C1)-(C3), then $I_f$ satisfies (DAC) with $A$ and $N$ if and only if there exists an automorphism $\varphi$ on [0,1] such that $af(b)\leq f(\varphi^{-1}(\varphi(a)+1-\varphi(\phi(N(b)))\wedge1))$ holds for any $a,b\in[0,1]$.
Especially, $I_f$ satisfies (A4) with $A$ and $N_I$ if and only if there exists an automorphism $\varphi$ on [0,1] such that
$af(b)\leq f(\varphi^{-1}(1+\varphi(a)-\varphi(\phi(\varphi^{-1}(1-\varphi(b))))\wedge1))$  holds for any $a,b\in[0,1]$.
\\{\bf Proof.} This proof is similar to that of Theorem 5.15.
\\{\bf Theorem 5.17} Let $I_g$ be a $g$-implication. If the continuous conjunctor $A$ fulfills (C1)-(C3), then $I_g$ satisfies (DAC) with $A$ and $N$ if and only if there exists an automorphism $\varphi$ on [0,1] such that $\frac{g(b)}{a}\leq g(\varphi^{-1}(\varphi(a)+1-\varphi(\phi(N(b)))\wedge1))$ holds for any $a,b\in[0,1]$.
Especially, $I_g$ satisfies (A4) with $A$ and $N_I$ if and only if there exists an automorphism $\varphi$ on [0,1] such that
$\frac{g(b)}{a}\leq g(\varphi^{-1}(1+\varphi(a)-\varphi(\phi(\varphi^{-1}(1-\varphi(b))))\wedge1))$  holds for any $a,b\in[0,1]$.
\\{\bf Proof.} It can be proved similarly to Theorem 5.15.
\\{\bf Theorem 5.18} Let $I_C$ be an probabilistic implication generated by the copula $C$ and the continuous conjunctor $A$ fulfill (C1)-(C3). $I_C$ satisfies (DAC) with $A$ and $N$ if and only if $C(a,b)\leq a\varphi^{-1}(\varphi(a)+1-\varphi(\phi(N(b)))\wedge1)$ holds, where $\varphi$ is an automorphism on [0,1].
Especially, $I_C$ satisfies (A4) with $A$ and $N_I$ if and only if
$C(a,b)\leq a\varphi^{-1}(1+\varphi(a)-\varphi(\phi(\varphi^{-1}(1-\varphi(b))))\wedge1)$.
\\{\bf Proof.} It can be proved similarly to Theorem 5.15.
\\{\bf Theorem 5.19} Let $\widetilde{I}_C$ be a probabilistic S-implication generated by the copula $C$. If the continuous conjunctor $A$ fulfills (C1)-(C3), then $\widetilde{I}_C$ satisfies (DAC) with $A$ and $N$ if and only if there exists an automorphism $\varphi$ on [0,1] such that $C(a,b)\leq \varphi^{-1}(\varphi(a)+1-\varphi(\phi(N(b)))\wedge1)+a-1$ holds for any $a,b\in[0,1]$.
Especially, $(\widetilde{I}_C,A)$ satisfies (A4) with $A$ and $N_I$ if and only if there exists an automorphism $\varphi$ on [0,1] such that
$C(a,b)\leq \varphi^{-1}(1+\varphi(a)-\varphi(\phi(\varphi^{-1}(1-\varphi(b))))\wedge1)+a-1$  holds for any $a,b\in[0,1]$.
\\{\bf Proof.} This proof is similar to that of Theorem 5.15.
\\{\bf Theorem 5.20} Let $I^T$ be $T$-power implication.

i. $I^{T_M}$ satisfies (DAC) with $A$ and $N$ if and only if $A$ has a right neutral element 1;

ii. If $T$ is a continuous Archimedean t-norm and the continuous conjunctor $A$ fulfills (C1)-(C3), then $I^T$ satisfies (DAC) with $A$ and $N$ if and only if there exists an automorphism $\varphi$ on [0,1] such that $\frac{t(a)}{t(b)}\leq \varphi^{-1}(\varphi(a)+1-\varphi(\phi(N(b)))\wedge1)$, where $t$ is the additive generator of $T$.
Especially, $I^T$ satisfies (A4) with $A$ and $N_I$ if and only if there exists an automorphism $\varphi$ on [0,1] such that
$\frac{t(a)}{t(b)}\leq \varphi^{-1}(1+\varphi(a)-\varphi(\phi(\varphi^{-1}(1-\varphi(b))))\wedge1)$  holds for any $a,b\in[0,1]$.
\\{\bf Proof.} This proof is similar to that of Theorem 5.15.
\section{Fuzzy reasoning with fuzzy implication satisfying (A4) ((A5)) and examples}
Let us firstly conclude our methods to discuss the satisfaction of (A4) or (A5) for fuzzy implications in this section. Further, we will propose an ACRI method employed the fuzzy implications which fulfill (A4) or (A5).
\subsection{Fuzzy reasoning with fuzzy implication satisfying (A4) or (A5)}
In order to investigate the the satisfaction of (A4) and (A5) for fuzzy implication, we use the following methods and processes:

(1) The following relationships are firstly revealed. That is,

i. $I$ satisfies (DAC) $\xrightarrow[]{\textmd{Proposition\ 3.3}}$ $I$ satisfies (A4);

ii. $I$ satisfies (AC) $\xrightarrow[]{\textmd{Proposition\ 3.9}}$ $I$ satisfies (A5).

Based on the relationship, we then discuss

(2) The satisfaction of (A5) for fuzzy implications.

i. For a given fuzzy implication $I$, $A_I$ can be constructed and $(I,A_I)$ satisfies (AC) (see Lemma 4.1);

ii. By Lemma 4.1 and Corollary 4.2, the condition for well-known fuzzy implications fulfilling (A5) can be obtained (see Theorems 4.3,4.6-4.9, 4.12, 4.14 and 4.15).

Parallelly, we finally study

(3) The satisfaction of (A4) for fuzzy implication.

i. If $I$ fulfills (CP(N)), then (A5) $\Longrightarrow$ (A4) (see Theorem 5.10). Therefore, it is sufficient to consider the condition for $I$ which satisfies (CP(N)) and (A5).

ii. If $I$ does not fulfill (CP(N)), the following two tactics are utilized:

(i) $I$ is respectively modified as $I^{lc}_{I,N}$ and $I^{uc}_{I,N}$ which satisfy (CP(N)) (see Theorem 5.11);

(ii) Demand that the aggregation function $A$ meets (C1)-(C3) (see Theorem 5.12). Especially, the conditions for well-known fuzzy implications satisfying (DAC) can be obtained (see Theorems 5.13-5.20).

It is well known that the fuzzy implication $I$ is usually  used to interpret the Promise 1 in FMP and FMT problems. According to the conclusion above, we thus propose the ACRI method as follows. The aggregation function $A$ would rather be selected such that $(A,I)$ fulfills (A5) in the ACRI method for FMP problem. This case implies that the conclusion $B'$ in FMP problem accords with the following intuition of human being: The nearer are between promise 2 and the antecedent of promise 1, the nearer are between conclusion and the consequent of promise 1. Especially, the conclusion is equal to the consequent of promise 1 if promise 2 is the antecedent of promise 1. For instance, if the R-implication $I_T$ generated by the left continuous t-norm $T$ is chosen to interpret the Promise 1 in FMP problem, we should select an aggregation function $A$ with a left neutral element 1 such that $A\leq T$ to compose the Promise 1 with Promise 2 in the ACRI method for FMP problem according to Theorem 4.3.

Similarly, we would rather choose the aggregation function $A$ such that $(A,I)$ fulfills (A4) in the ACRI method for FMT problem. In this case, the conclusion $D'$ in FMT problem conforms to the following intuition of human being: The farther are between promise 2 and the consequent of promise 1, the farther are between conclusion and the antecedent of promise 1. Especially, the conclusion is equal to the negation of consequent of promise 1 if promise 2 is the negation of antecedent of promise 1. For example, if the $(A,N)$-implication come from  a continuous
t-conorm $S$ with the ordinal sum structure  and a strong negation $N$ is chosen to interpret the Promise 1 in FMT problem,  the aggregation function $A$ should be selected such that $A\leq A_{I_{S,N}}$ to compose the Promise 1 with Promise 2 in ACRI method for FMT problem by Theorem 5.10. Especially, for a given aggregation function $A$ in the ACRI method for FMT problem, if we know partial information of $A$, that is, $A$ fulfills (C1)-(C3), then we can select a suitable fuzzy implication to interpret the Promise 1 in FMT problem by Theorems 5.13-5.20.
\subsection{Examples}
In this subsection we will present two examples to illustrate the results obtained in the previous sections.
\\{\bf Example 6.1} To automatically diagnose whether the patients  have illness, we need to construct a flexible classification to divide entities which denote the patients possess the attributions into three classes, patient has a
diagnosis (Y), probably has (further evaluation is advisable) (M), and
does not have (N). According to the knowledge of domain experts,  a fuzzy medical classification is formalized by two attributes as follows:

IF attribute 1 is\ $high$\ AND attribute 2 is\ $high$,\ THEN  \ Y,

IF attribute 1 is\ $medium$\ AND attribute 2 is\ $high$,\ THEN  \ Y,

IF attribute 1 is\ $low$\ AND attribute 2 is\ $high$,\ THEN  \ M,

IF attribute 1 is\ $high$\ AND attribute 2 is\ $low$,\ THEN  \ M,

IF attribute 1 is\ $medium$\ AND attribute 2 is\ $low$,\ THEN  \ N,

IF attribute 1 is\ $low$\ AND attribute 2 is\ $low$,\ THEN  \ N.

Where attributes 1 and 2 denote a patient has exhibited symptoms. Here, we use the following linguistic variables for attribute 1: high, medium and low. And the attribute 2 is characterized with the following linguistic labels: high and low. The membership functions of linguistic variables for attributes 1 and 2 are shown in Fig.1.
\begin{center}
\includegraphics[width=7cm, height=4cm]{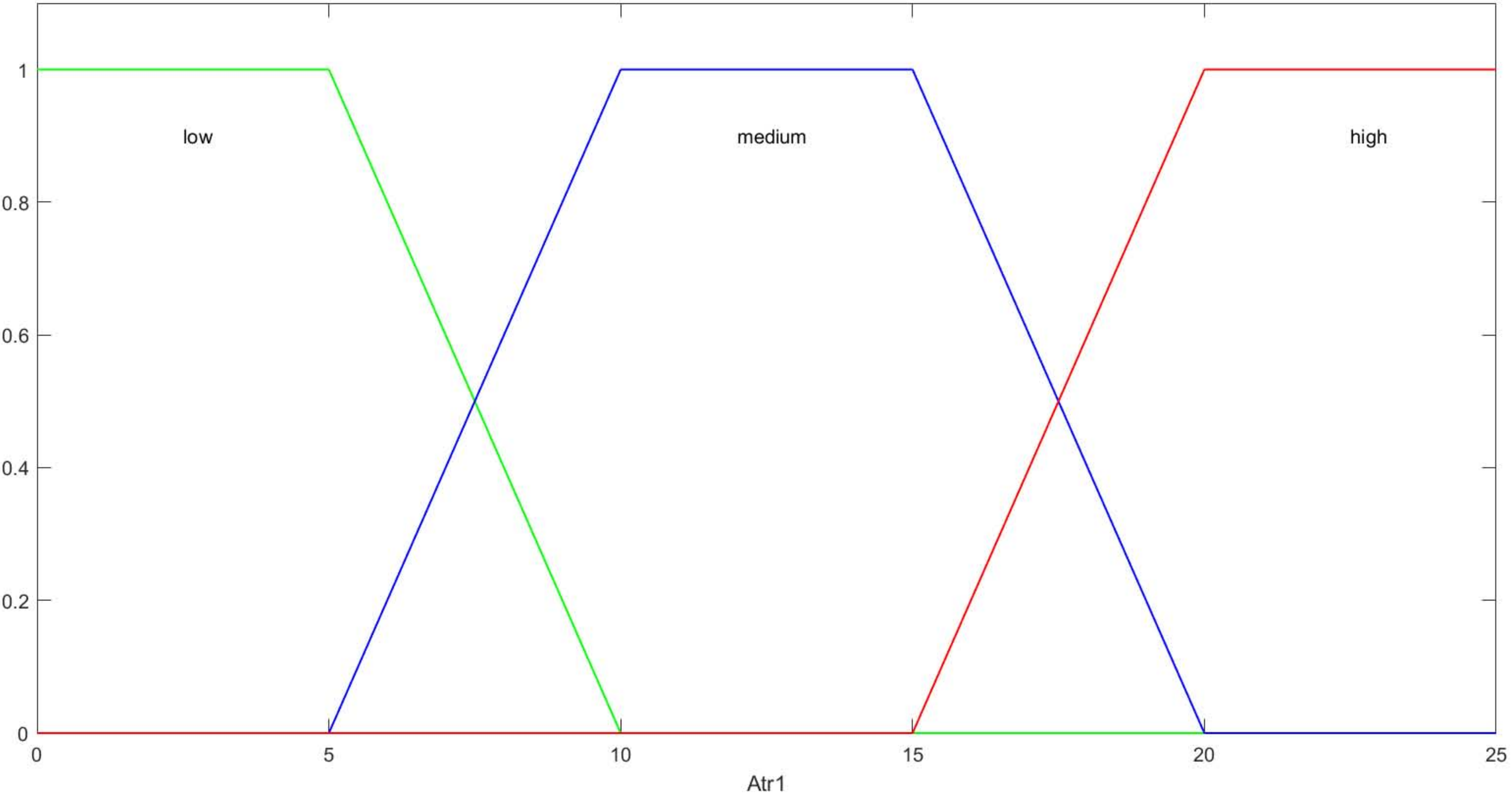}\vspace{1mm}
\includegraphics[width=7cm, height=4cm]{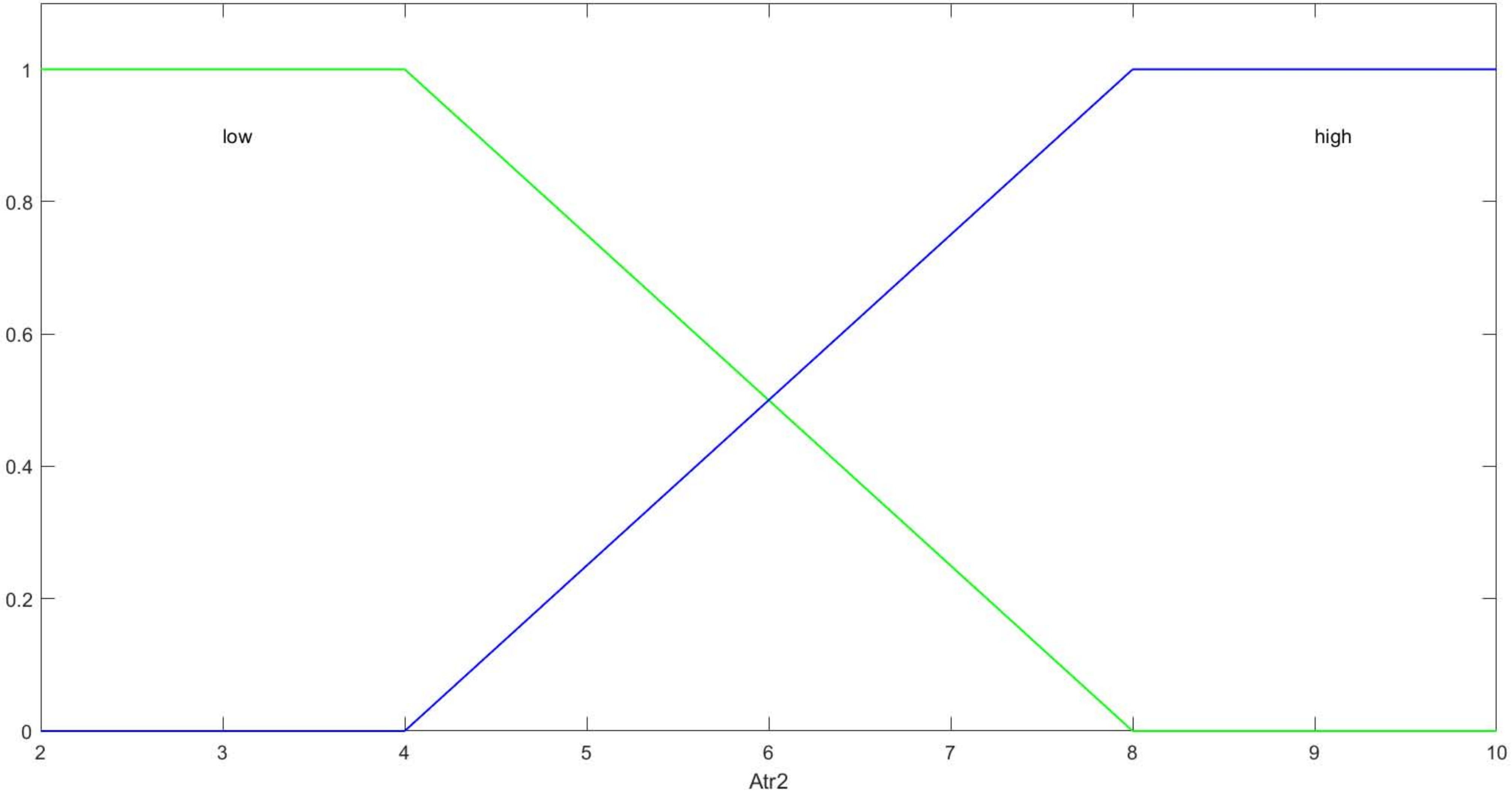}
\vspace*{-.3cm} $$\mbox{{\small\ {\bf\small Figure 1}
\quad{\bf\small The membership functions of linguistic variables for attributes 1 and 2 \qquad}}}$$
\end{center}

 To obtain the conclusion from the following pattern\vspace{-1mm}
$$\textmd{Premise\ 1:}\ \textmd{IF}\ x\ \textmd{is }D\ \textmd{THEN}\ y\ \textmd{is}\ B\qquad\vspace{-2mm}$$
$$\textmd{Premise\ 2:}\ \ x\ \textmd{is }D'\qquad\qquad\qquad\qquad\quad\vspace{-1mm}$$
$$D'\ \textmd{differs\ from D,\ but} \ D'\ \textmd{is\ near\ from}\ D\vspace{-3mm}$$
\begin{picture}(1,2)
 \qquad\qquad\qquad\qquad\line(1,0){220}
\end{picture}\vspace{-1mm}
$$\qquad\textmd{Conclusion:}\ y\ \textmd{is}\ B'\ \textmd{with}\ B' \ \textmd{is\ not\ far\ from}\ B,$$
 the aggregation function $A(x,y)=(x+y-1)\vee 0$ is employed to interpret the word ``and" while the {\L}ukasiewicz implication $I_{LK}(x,y)=(1-x+y)\wedge 1$ is used to interpret the fuzzy rule. By the discussion above,  they fulfill (A5).  For the entity $E_1=(11, 3)$, we use triangular fuzzifier to translate $E_1$ into a fuzzy input which attribute 1 is medium and attribute 2 is low. Obviously, $E_1$ satisfies the antecedent of the fifth fuzzy rules. With our proposed method, the conclusion should be equal to the consequent of the fifth fuzzy rules. Therefore,  we have the fact that $E_1$ belongs to the class ``N". Similarly, the entity $E_2=(20, 2)$ belongs to the class ``M". And the entity $E_3=(22, 8)$ belongs to the class ``Y". These entities belonging to the classes are shown in Figure 2.
\begin{center}
\includegraphics[width=8cm, height=6cm]{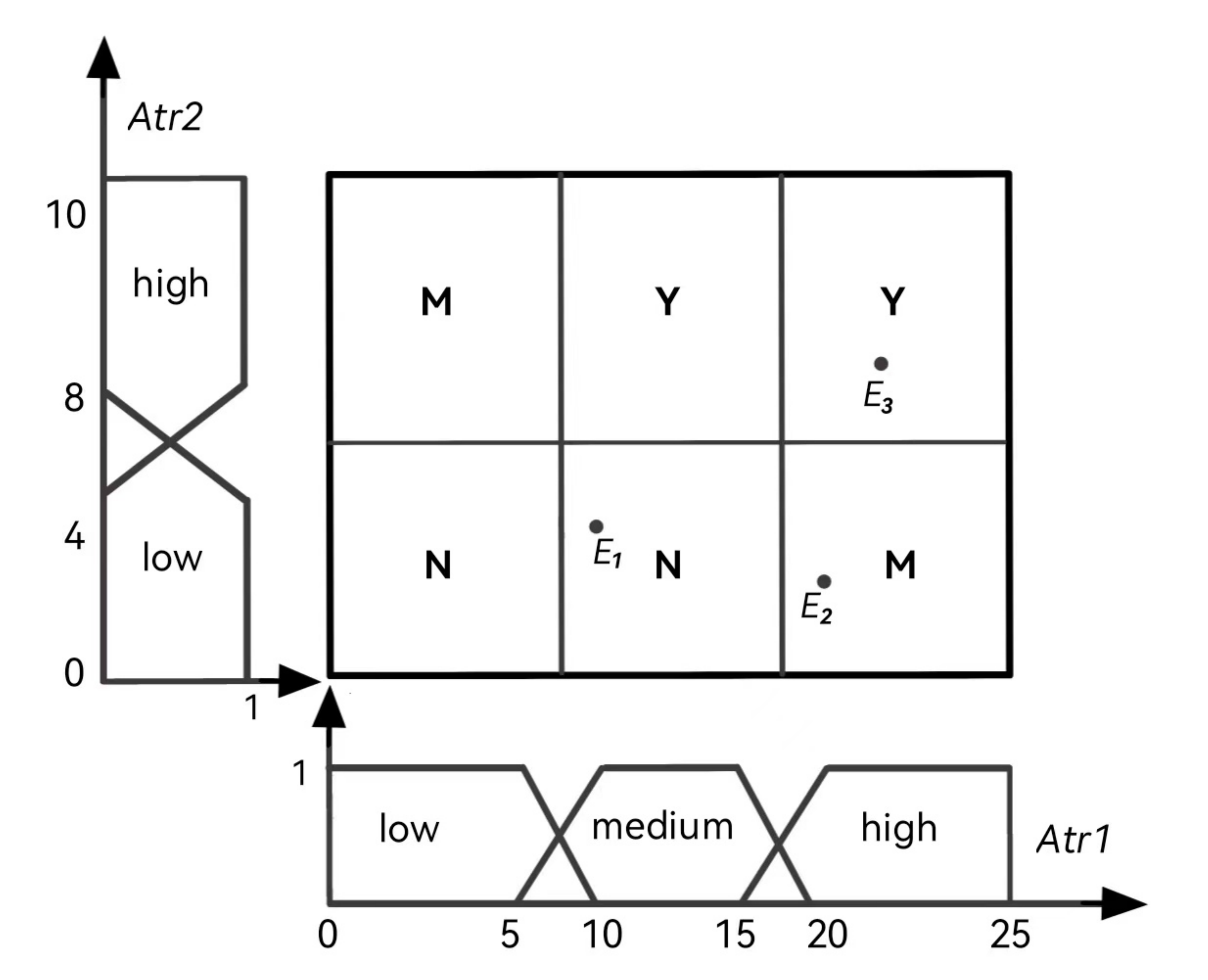}
 \vspace*{-.3cm} $$\mbox{{\small\ {\bf\small Figure 2}
\quad{\bf\small The classes of entities with attributes 1 and 2 \qquad}}}$$
\end{center}
{\bf Example 6.2}  Let $D_1=high=[1,0.1,0,0.05]$ be a fuzzy set on $U_1=\{x_{11},x_{12},x_{13},x_{14}\}$ denoted the sugar of apple, $D_2=red=[0,0.9,0.04,0,0]$ a fuzzy set on $U_2=\{x_{21},x_{22},x_{23},x_{24},x_{25}\}$ expressed the colour of apple and $B=ride=[0.3,0.2,0.4]$  a fuzzy on $V=\{y_1,y_2,y_3\}$ described the maturity of apple. In order to deduce the apple whether is ride, we suppose to consider the two-input-one-output fuzzy system including the following single fuzzy IF-THEN rule

 \qquad\qquad\qquad IF \ $x_1$\ is\ $D_{1}$ and \ $x_2$\ is \ $D_{2}$ THEN  \ $y$\ is \ $B$.

Assume that the Reichenbach implication
$I_{RC}(x,y)=1-x+xy$ is used to interpret the IF...THEN rule. Similarly, in order to achieve the scheme of fuzzy reasoning mentioned in Example 6.1, we construct the aggregation\vspace{1mm} function $A_{I_{RC}}$ as $A_{I_{RC}}(x,y)=\left\{\begin{array}{ll}
                   1-\frac{1-y}{x}& x\neq 0\ \textmd{and}\ x+y\geq 1\\
                   0&\textmd{otherwise}
                 \end{array}\right.$ according to Lemma 4.1. By Theorem 4.3, we have the fact that $I_{RC}$ satisfies (A5) with $A_{I_{RC}}$. Let $D'_1=$\vspace{1mm} $[1,0,0,0]$ and $D'_2=[0,1,0,0,0]$ be the fuzzy single input. With the ACRI method involved $I_{RC}$ and $A_{I_{RC}}$,  the conclusion $B'$ can be calculated as $B'=[0.37,0.28,0.46]$.

                 We can assert that $B'$ and $B$ is near because the fuzzy input $D'=D'_1\times D'_2$ is close to the antecedent of fuzzy rule $D=D_1\times D_2$. Indeed, according to the distance measures between two fuzzy sets in \cite{Zwick}, the distance between $D$ and $D'$ is defined as $d(D',D)=\left(\sum\limits_{i=1}^n|D(x_i)-D'(x_i)|^p\right)^{\frac{1}{p}}$\vspace{1mm} with $p\geq 1$. Let $p=2$.\vspace{1mm} We can obtain $d(D'_1\times D'_2, D_1\times D_2)\doteq0.108$ and $d(B,B')\doteq0.122$. Especially,  the output $B'=B$ if $D'=D'_1\times D'_2=D_1\times D_2=D$. However, if we chose the greatest disjunctor $D_\top$ to interpret the word ``and" in FMP problem. It is easy to see that $I_{RC}$ does not fulfills (A5) with $D_\top$. And then the conclusion $B''$ of FMP problem with the ACRI method involved $I_{RC}$ and $D_\top$ is $B''=[1,1,1]$. Obviously, $d(B,B'')\doteq1.22>0.122\doteq d(B,B')$.
\\{\bf Remark 7.} For convenience to show our proposed method, Example 6.2 only consider the case that the fuzzy sets defined on the  discrete universes. It is not difficult to see that the conclusions are near if the fuzzy input and the antecedent of fuzzy rule are not far in the case that the fuzzy sets define on the  continuous universes. Especially, that the conclusions are $B$ ($D^C)$ if $D'=D$ ($D'=B^C$). This means that our proposed method possesses the adaptability for any fuzzy set involved in FMP and FMT problems. Therefore, it can be utilized to make the valid uncertain reasoning in practice.
\subsection{Advantages and limitations of our proposed method}

From the above example, if the fuzzy implication and aggregation function satisfy (A5), it is clear that the output of FMP  problem accords with the following intuition of human being: The nearer are between promise 2 and the antecedent of promise 1, the nearer are between conclusion and the consequent of promise 1. Especially, the conclusion is equal to the consequent of promise 1 if promise 2 is the antecedent of promise 1. And then
the advantages of our proposed inference method can be summarized as follows.
\begin{itemize}
  \item Satisfaction of classical modus ponens. For a fuzzy implication which is used to interpret the fuzzy rule, we can construct a corresponding aggregation function such that the output of FMP problem is $B'=B$ when the input is equal to the antecedent of fuzzy rule.
  \item Satisfaction of classical modus tollens. If a given fuzzy implication is used to interpret the fuzzy rule, a corresponding aggregation function can be obtained such that the output of FMT problem  is $B'=D^C$ when the input is the negation
of antecedent of promise 1. Especially, for a given aggregation function $A$ which is used to interpret the word ``and", if $A$ fulfills (C1)-(C3), then we can select a suitable fuzzy implication to interpret the Promise 1 in FMT problem, and then the the output of FMT problem is $B'=D^C$ when the input is the negation
of antecedent of promise 1.
\end{itemize}

Although the classical modus ponens has been consider in some methods to solve the FMP problem\cite{Liz,Liw,Magrez,Pei,Raha,Wang1,Zhou}, these methods cannot ensure that the classical modus tollens property holds. By contrast, our proposed method considers not only the classical modus ponens property but also the classical modus tollens property. However, our proposed method also has the following limitations.
\begin{itemize}
  \item Deficiency of a logical foundation. Considering our proposed method based on ACRI method, it still possesses the deficiency of CRI method as pointed out by some researchers\cite{Baldwin,Mizumoto,Turksen,Wang1,Zhou}.
  \item Failed to ensure the classical modus ponens and modus tollens hold together. Unlike the methods in Ref.\cite{Gera,Trillas}, our proposed method cannot find an aggregation function such that the classical modus ponens and  modus tollens hold together for the well-known fuzzy implications. Notice that Gera used a sigmoid-like function to construct fuzzy set in order to ensure the reasoning fulfills MP and MT properties together. Perhaps, this is a direction of our research in future.
\end{itemize}
With comparing the different fuzzy reasoning
methods with our proposed method, we finally list together these fuzzy reasoning methods in Table 4.
 $$\mbox{\bf{\small Table\ 4 \ Comparion of several fuzzy resoning methods}}\vspace{-3mm}$$
\begin{center}
\begin{tabular}{lccccc}
 \toprule[1pt]
&Expressing the word``AND"&Interpreting fuzzy rule&MP & MT & Logic foundation\vspace{1mm}\\
\midrule[0.75pt]
Baldwin's method &t-norm&{\L}ukasiewicz Implication&No&No&Yes\\
Mizumoto's method &{\L}ukasiewicz t-norm&{\L}ukasiewicz Implication&Yes&Yes&Yes\\
Turksen's method  & Matching function&Modification function&Yes&No&Yes\\
Wang's method&Left-continuous t-norm&R-implication&Yes&No&Yes\\
Zhou's method  &Left-continuous t-norm&R-implication&Yes&No&Yes\\
Gera's method &t-norm&Fuzzy implication&Yes&Yes&No\\
Trillas' method &Continuous t-norm&R (S)-implications&Yes&Yes&No\\
Our proposed method&Aggregation functions&More fuzzy implications&Yes&Yes&No\\
   \bottomrule[1pt]
\end{tabular}
\end{center}\vspace{1mm}
\section{Conclusions}
 Aggregation functions and fuzzy implications play an important role in fuzzy inference. In order to enhance the effectiveness of ACRI method, we then have studied the MP and MT properties of the ACRI method with well-known fuzzy implications using the axioms (A4) and (A5) in detail. Concretely, we have

(1)  Analyzed the properties of aggregation functions and fuzzy implications involved satisfying (A4) or (A5) in the ACRI method;

(2) Constructed the aggregation function for well-known fuzzy implications satisfying (A5) in the ACRI method;

(3) Given the conditions for well-known fuzzy implications satisfying (A4) with a strong negation in the ACRI method;

(4) Shown an ACRI fuzzy inference method involved well-known fuzzy implications and aggregation functions which fulfill (A4) or (A5).

These results contribute to improve the effectiveness of the ACRI method. In the future, we
will investigate the validity of the ACRI method using the axioms (A6) and (A7). Considering that some linguistic modifiers are involved in (A6) and (A7), we will extend them as follows:

(A6$'$) $B'=m(B)$ if $D'=m(D)$;

(A7$'$) $B'=B$ if $D'=m(D)$, where $m$ is a fuzzy modifier.

We also will study the MP and MT properties of ACRI method for interval-valued fuzzy sets, fuzzy soft sets and so on. And will apply them in real-life decision-making.

Moreover, considering there exists still some deficiencies in ACRI method as pointed out by
some researchers \cite{Baldwin,Mizumoto,Turksen,Wang,Zhou}, we will study the MP and MT properties for similarity-based approximate reasoning (SBR) method, triple implication principle (TIP) method and  quintuple implication principle (QIP).
\section{Acknowledgement}
The authors would like to thank the anonymous referees and the Editor-in-Chief for their valuable comments. This work was supported by the National Natural Science Foundation of China (Grant No. 61673352).
\section{Compliance with ethical standards}
\noindent{\bf Conflict of interest} Author declares that he has no conflict of interest.
\vspace{1mm}
\\{\bf Human and animal rights} This article does not contain any studies
with human participants or animals performed by the authors.
\vspace{1mm}
\\{\bf Data availability statement} This article has no associated data.

\end{document}